\documentclass[12pt]{amsart}
\usepackage{amsfonts}
\usepackage{amssymb}
\newtheorem{thm}{Theorem}[section]

\newtheorem{lem}[thm]{Lemma}
\newtheorem{prop}[thm]{Proposition}
\theoremstyle{definition}
\newtheorem{defn}[thm]{Definition}
\theoremstyle{remark}

\newcommand{\R}{\mathbb R}
\newcommand{\C}{\mathbb C}
\newcommand{\D}{\mathbb D}
\newcommand{\T}{\mathbb T}

\newcommand{\X}{\mathfrak X}

\begin{document}
\title[K\"ahler metrics whose potential is a function of
the time-like distance] {K\"AHLER METRICS GENERATED BY FUNCTIONS
OF THE
TIME-LIKE DISTANCE IN THE FLAT K\"AHLER-LORENTZ SPACE}%
\author{G. Ganchev and V. Mihova}%
\address{Bulgarian Academy of Sciences, Institute of Mathematics
and Informatics,Acad. G. Bonchev Str. bl. 8, 1113 Sofia, Bulgaria}%
\email{ganchev@math.bas.bg}%
\address{Faculty of Mathematics and Informatics, University of Sofia,
J. Bouchier Str. 5, (1164) Sofia, Bulgaria}
\email{mihova@fmi.uni-sofia.bg}

\begin{abstract}
We prove that every K\"ahler metric, whose potential is a function
of the time-like distance in the flat K\"ahler-Lorentz space, is
of quasi-constant holomorphic sectional curvatures, satisfying
certain conditions. This gives a local classification of the
K\"ahler manifolds with the above mentioned metrics. New examples
of Sasakian space forms are obtained as real hypersurfaces of a
K\"ahler space form with special invariant distribution. We
introduce three types of even dimensional rotational hypersurfaces
in flat spaces and endow them with locally conformal K\"ahler
structures. We prove that these rotational hypersurfaces carry
K\"ahler metrics of quasi-constant holomorphic sectional
curvatures satisfying some conditions, corresponding to the type
of the hypersurfaces. The meridians of those rotational
hypersurfaces, whose K\"ahler metrics are Bochner-K\"ahler
(especially of constant holomorphic sectional curvatures) are also
described.
\end{abstract}
\maketitle

\section{Introduction}

In \cite {GM2} we have given a complete description of the
curvature tensor and curvature properties of the K\"ahler metrics
$g = \partial\bar\partial f(r^2)$, where $r^2$ is the distance
function with respect to the origin in ${\C}^n$ and the real
$\mathcal{C}^{\infty}$-function $f(r^2)$ satisfies the conditions
$$ f'(r^2) > 0, \quad f'(r^2) + r^2f''(r^2) > 0. $$

Bochner-K\"ahler metrics of the type  $\partial\bar\partial
f(r^2)$ have been studied in \cite {TL}. The completeness of these
metrics has been discussed in \cite {B}.

We have introduced the notion of a K\"ahler manifold $(M, g, J, D)
\,(\dim M = 2n \geqq 6)$ with $J$-invariant $B_0$-distribution
$D\, (\dim D = 2(n-1)).$ Any $B_0$-distribution generates a
function $k > 0$ on $M$. If $D^{\perp}$ is the distribution,
orthogonal to $D$, then every holomorphic section $E(p),\, p \in
M,$ determines a geometric angle $\vartheta = \angle (E(p),
D^{\perp} (p)).$

A K\"ahler manifold $(M, g, J, D)$ is of quasi-constant
holomorphic sectional curvatures if its holomorphic sectional
curvatures only depend on the point $p$ and the angle $\vartheta$.

If $(M, g, J, D)$ is a K\"ahler manifold of quasi-constant
holomorphic sectional curvatures, then the distribution
$D(p),\,p\in M$ is of pointwise constant holomorphic sectional
curvature $a(p)$ and the function $ a+k^2$ divides the class of
these manifolds into three subclasses according to
$$a+k^2>0, \quad a+k^2=0, \quad a+k^2<0.$$

In \cite {GM2} we have shown that the flat K\"ahler manifold
${\C}^n$ carries a canonical $B_0$-distribution and proved the
following characterization of the family of K\"ahler metrics $g =
\partial\bar\partial f(r^2)$: \vskip 1mm {\it Any K\"ahler metric
$g = \partial\bar\partial f(r^2)$ is of quasi-constant holomorphic
sectional curvatures with $ a + k^2 > 0$.

Conversely, any K\"ahler manifold $M \, (\dim M=2n \geq 6)$ of
quasi-constant holomorphic sectional curvatures with
$B_0$-distribution and $a+k^2>0$ is locally equivalent to
$({\C}^n,g,J_0)$ with the canonical $B_0$-distribution and
$g=\partial\bar\partial f(r^2)$.} \vskip 2mm In this paper we
solve the problem of describing the curvature properties of the
K\"ahler metrics generated by potential functions $f(-r^2)$,
$-r^2$ being the time-like distance function from the origin in
the flat K\"ahler-Lorentz space.

Let $({\C}^n, h', J_0)$ be the flat K\"ahler-Lorentz space with
the canonical complex structure $J_0$ and flat K\"ahler metric
$h'$ of signature $(2(n-1),2)$.

In Proposition \ref{P:3.2} we prove that if $f(-r^2),\, -r^2 < 0,$
is a real $\mathcal{C}^{\infty}$-function satisfying the
conditions
$$ f'(-r^2) > 0, \quad f'(-r^2) -r^2 f''(-r^2) < 0, $$
then $ g = \partial\bar\partial f(-r^2) $ is a positive definite
K\"ahler metric in the time-like domain ${\T}_1^{n-1}= \{
\textbf{Z} \in {\C}^n :\, h'(\textbf{Z}, \textbf{Z}) < 0 \}.$

In Section 4 we prove the basic Theorem \ref{T:4.1}, which gives a
complete curvature description of the family of K\"ahler metrics $
g = \partial\bar\partial f(-r^2) $: \vskip 1mm {\it Any K\"ahler
metric $g = \partial\bar\partial f(-r^2)$ is of quasi-constant
holomorphic sectional curvatures with $a+k^2<0$.

Conversely, every K\"ahler manifold $M \,(\dim M= 2n \geq 6)$ of
quasi-constant holomorphic sectional curvatures with
$B_0$-distribution and $a+k^2<0$ is locally equivalent to
$({\T}^{n-1}_1,g,J_0)$ with the canonical $B_0$-distribution and
$g=\partial\bar\partial f(-r^2)$.} \vskip 1mm In Section 5 we
clear up the geometric meaning of the function $a+k^2$ in a
K\"ahler manifold $(M,g,J,D)$ of quasi-constant holomorphic
sectional curvatures. We show that $(M,g,J,D)$ is a one-parameter
family of $\alpha$-Sasakian space forms $Q^{2n-1}(s), \; s \in I$
with $\alpha =\displaystyle{\frac{k}{2}}$ and prove in Proposition
\ref{P:5.2} that $sign (a+k^2)$ determines the type of the
corresponding $Q^{2n-1}(s)$.

As a consequence of Theorem \ref{T:4.1} we obtain examples of
K\"ahler space forms in ${\T}^{n-1}_1$ with $B_0$-distribution and
$a+k^2<0$. Especially the metric $ g = - 2\partial\bar\partial \ln
\,(r^2-1), \quad -r^2  < -1,$ is of constant holomorphic sectional
curvature  $-1$. Considering the unit "disc" $({\D}^{n-1}_1(1) :
h'(\textbf{Z}, \textbf{Z})<-1)$ we show that any hypershere
$H^{2n-1}_1(O,r), \; r>1$ in $({\D}^{n-1}_1,g,J_0)$ carries a
natural structure of an $\alpha$-Sasakian space form with
$\displaystyle{\alpha=\frac{1}{2r}}$ and constant
$\varphi$-holomorphic sectional curvatures $c$, so that
$c+3\alpha^2 <0$ (cf {\cite{T}).

In Section 6 we consider three types of rotational hypersurfaces
$M$ in ${\C}^n\times{\R}$ with axis of revolution $l={\R}$:

I type: the parallels $S^{2n-1}$ are the usual hyperspheres in the
complex Euclidean space $({\C}^n,g',J_0)$ and the axis of
revolution ${\R}$ is endowed with positive definite inner product;
the meridians are curves in Euclidean plane.

II type: the parallels $S^{2n-1}$ are the usual hyperspheres in
the complex Euclidean space $({\C}^n,g',J_0)$ and the axis ${\R}$
is endowed with negative definite inner product; the meridians are
space-like curves in hyperbolic plane.

III type: the parallels $H^{2n}_1$ are hyperspheres in the flat
time-like domain $({\T}^{n-1}_1,h',J_0)$ and the axis ${\R}$ is
endowed with positive definite inner product; the meridians are
time-like curves in hyperbolic plane.

In Subsection 6.1 we recall that the hypersurfaces of type I carry
a natural K\"ahler structure of quasi-constant holomorphic
sectional curvatures with functions $a>0, \; a+k^2>0$. In
Proposition \ref{P:6.3} we obtain the meridians of the rotational
hypersurfaces of type I, whose K\"ahler metric is
Bochner-K\"ahler.

In Subsection 6.2 we introduce a K\"ahler structure on rotational
hypersurfaces of type II and prove in Theorem \ref{T:6.2} that
this K\"ahler structure is of quasi-constant holomorphic sectional
curvatures with functions $a<0, \; a+k^2>0$. We find the meridians
of the rotational hypersurfaces of type II, whose K\"ahler metric
is Bochner-K\"ahler (Proposition \ref{P:6.6}) or of constant
holomorphic sectional curvatures (Proposition \ref{P:6.5})

In Subsection 6.3 we introduce a K\"ahler structure on the
rotational hypersurfaces of type III and prove in Theorem
\ref{T:6.3} that this K\"ahler structure is of quasi-constant
holomorphic sectional curvatures with functions $a<0, \; a+k^2<0$.
We find the meridians of those rotational hypersurfaces of type
III, whose K\"ahler metric is Bochner-K\"ahler (Proposition
\ref{P:6.8}) or is of constant holomorphic sectional curvatures
(Proposition \ref{P:6.7}). \vskip 2mm
\section{Preliminaries}
\vskip 2mm

In this section we give some basic notions and formulas for
K\"ahler manifolds with $B_0$-distribution \cite {GM2} we need
further.

Let $(M,g,J,D)$ be a 2n-dimensional K\"ahler manifold with metric
$g$, complex structure $J$ and $J$-invariant distribution $D$ of
codimension 2. The Lie algebra of all $C^{\infty}$ vector fields
on $M$ will be denoted by ${\X}M$ and $T_pM $ will stand for the
tangent space to $M$ at any point $p \in M$. In the presence of
the distribution $D$ the structure of any tangent space is $T_pM =
D(p) \oplus D^{\perp}(p)$, where $D^{\perp}(p)$ is the
2-dimensional $J$-invariant orthogonal complement to the space
$D(p)$. This means that the structural group of the manifolds
under consideration is the subgroup $U(n-1)\times U(1)$ of $U(n)$.

In the local treatment of these manifolds $D^{\perp} = span \{\xi,
J\xi\}$ for some unit vector field $\xi$. The 1-forms,
corresponding to $\xi$ and $J\xi$, respectively, are
$$\eta(X) = g(\xi, X), \quad \tilde \eta (X) = g(J\xi, X)
= -\eta(JX); \quad X \in {\X}M.$$
Then the distribution $D$ is determined by the conditions
$$D(p) = \{X \in T_pM \, \vert \, \eta(X) = \tilde \eta(X)
= 0 \}, \quad p \in M.$$

The K\"ahler form $\Omega$ of the structure $(g,J)$ is given by
$\Omega(X,Y) = g(JX,Y)$, $X,Y \in {\X}M$.

Let $\nabla$ be the Levi-Civita connection of the metric $g$.

A J-invariant distribution
$D,\, (D^{\perp} = span \{\xi, J\xi \})$ is said to be a
$B_0$-distribution \cite {GM2} if
$$\begin{array}{l}
i) \displaystyle{ \quad \nabla _{x_0} \xi =
\frac{k}{2}\,x_0, \; k\neq 0;}\quad x_0 \in D;\\
[2mm]
ii) \quad \nabla_{J\xi}\xi=-p^*J\xi;\\
[2mm]
iii)\quad \nabla _{\xi} \xi = 0. \\
[2mm]
\end{array}$$

The above definition implies immediately the following equalities
\cite {GM2}
$$\nabla_X\xi=\frac{k}{2}\{X-\eta(X)\xi-\tilde\eta(X)J\xi\}
-p^*\tilde\eta(X)J\xi, \quad X\in {\X}M; \leqno (2.1)$$
$$dk = \xi(k)\,\eta, \quad p^* = -\frac{\xi(k)+k^2}{k};\leqno{(2.2)}$$

Any K\"ahler manifold $(M, g, J, D)$ with $J$-invariant
distribution $D$ carries the tensors
$$\begin{array}{ll}
4\pi (X,Y)Z := &g(Y,Z)X - g(X,Z)Y - 2g(JX,Y)JZ\\[2mm]
&+ g(JY,Z)JX - g(JX,Z)JY;
\end{array}\leqno{(2.3)}$$
$$\begin{array}{l}
\Phi_1(X,Y)Z:=\\
[2mm]
\displaystyle{\frac{1}{8}\{g(Y,Z)(\eta(X)\xi+\tilde\eta(X)J\xi)
-g(X,Z)(\eta(Y)\xi+\tilde\eta(Y)J\xi)}\\
[2mm]
+g(JY,Z)(\eta(X)J\xi-\tilde\eta(X)\xi)-g(JX,Z)(\eta(Y)J\xi
-\tilde\eta(Y)\xi)\\
[2mm]
-2g(JX,Y)(\eta(Z)J\xi-\tilde\eta(Z)\xi)\};\\
[2mm]
\Phi_2(X,Y)Z:=\\
[2mm]
\displaystyle{\frac{1}{8}\{(\eta(Y)\eta(Z)+
\tilde\eta(Y)\tilde\eta(Z))X
-(\eta(X)\eta(Z)+\tilde\eta(X)\tilde\eta(Z))Y}\\
[2mm]
+(\eta(Y)\tilde\eta(Z)-\tilde\eta(Y)\eta(Z))JX
-(\eta(X)\tilde\eta(Z)-\tilde\eta(X)\eta(Z))JY\\
[2mm]
-2(\eta(X)\tilde\eta(Y)-\tilde\eta(X)\eta(Y))JZ\};\\
[1mm]
\Phi := \Phi_1 + \Phi_2;
\end{array}\leqno{(2.4)}$$
$$\begin{array}{lll}
\Psi (X,Y)Z& := &\eta(Y)\eta(Z)\tilde\eta(X)J\xi
- \eta(X)\eta(Z)\tilde\eta(Y)J\xi\\
[2mm]
& &+ \eta(X)\tilde\eta(Y)\tilde\eta(Z)\xi -
\eta(Y)\tilde\eta(X)\tilde\eta(Z)\xi\\
[2mm]
&=&(\eta \wedge\tilde\eta)(X, Y)(\tilde\eta(Z)\xi -
\eta(Z)J\xi),
\end{array}\leqno{(2.5)}$$
$X,Y,Z \in {\X}M.$ These tensors are invariant under the action of
the structural group $U(n-1)\times U(1)$ \cite {TV}.

The Riemannian curvature tensor $R$ of the metric $g$ is given by
$$R(X,Y)Z = \nabla_X\nabla_YZ - \nabla_Y\nabla_XZ - \nabla_{[X,Y]}Z,$$
$$R(X,Y,Z,U) = g(R(X,Y)Z,U); \quad X,Y,Z,U \in {\X}M.$$

In \cite {GM2} we proved that a K\"ahler manifold $(M, g, J, D)\;
(\dim M = 2n \geq 4)$ with $J$-invariant distribution $D$ is of
quasi-constant holomorphic sectional curvatures if and only if
$$R = a\pi + b \Phi + c\Psi,$$
where $a$, $b$ and $c$ are functions on $M$, generated by the
structure $(g,J,\xi)$.

If $(M, g, J, D)\; (\dim M = 2n \geq 6)$ is a K\"ahler manifold of
quasi-constant holomorphic sectional curvatures, then the
following statements hold good \cite {GM2}:

(i)\quad\; \,If $D$ is a $B_0$-distribution, then
$$ \displaystyle { da = \frac{kb}{2}\eta.} \leqno (2.6)$$

(ii)\quad\;Under the condition $b \neq 0$ \,$D$ is a
$B_0$-distribution if and only if $D$ is non-involutive.

(iii) \quad If $b = 0$ and $D$ is non-involutive, then $c = 0$,
i.e. $M$ is a K\"ahler space form.

Finally we recall some basic facts related to $\alpha$-Sasakian
manifolds.

Let $Q^{2n-1}(g, \varphi, \tilde\xi, \tilde\eta)\, (n \geq 3)$ be
an almost contact Riemannian manifold, i.e.
$$\begin{array}{l}
\vspace{2mm} g(\varphi x, \varphi y)=g(x,
y)-\tilde\eta(x)\tilde\eta(y),
\; x, y \in {\X}Q^{2n-1},\\
\vspace{2mm}
{\varphi}^2 x=-x+\tilde\eta(x)\tilde\xi,\; x \in {\X}Q^{2n-1},\\
\vspace{2mm} \varphi \, \xi=0.
\end{array} \leqno (2.7)$$

If the structure $(g, \varphi, \tilde\xi, \tilde\eta)$ of an
almost contact Riemannian manifold $Q^{2n-1}$ satisfies the
conditions
$$\mathcal{D}_x\tilde\xi=\alpha \, \varphi x, \quad x \in {\X}Q^{2n-1},$$
$$(\mathcal{D}_x \varphi)(y)=\alpha \, (\tilde\eta(y)x-g(x, y)
\tilde\xi),\quad x, y \in {\X}Q^{2n-1},$$ where $\mathcal{D}$ is
the Levi-Civita connection of the metric $g$ and $\alpha = const$,
then $Q^{2n-1}$ is called an {\it $\alpha$-Sasakian manifold}
\cite {JV}.

If the constant $\alpha =1$, then $Q^{2n-1}$ is a Sasakian
manifold in the usual sense.

$\alpha$-Sasakian space forms are characterized as follows:
\begin{prop}\label{P:2.1}
$($\cite {O}, \cite {JV}$)$ An $\alpha$-Sasakian manifold
$(Q^{2n-1}, g, \varphi, \tilde\xi, \tilde\eta)\; (\dim
Q^{2n-1}\geq 5)$ is of constant $\varphi$-holomorphic sectional
curvatures $c$ if and only if
$$\begin{array}{ll}
K(x,y,z,u)=&\displaystyle{\frac{c+3{\alpha}^2}{4}}
[g(y,z)g(x,u)-g(x,z)g(y,u)]\\
[2mm] &+\displaystyle{\frac{c-{\alpha}^2}{4}}[g(\varphi
y,z)g(\varphi x,u)-
g(\varphi x,z)g(\varphi y,u)-2g(\varphi x,y)g(\varphi z,u)\\
[2mm]
&-g(y,z)\tilde\eta(x)\tilde\eta(u)-g(x,u)\tilde\eta(y)\tilde\eta(z)\\
[2mm]
&+g(x,z)\tilde\eta(y)\tilde\eta(u)+g(y,u)\tilde\eta(x)\tilde\eta(z)],
\quad x,y,z,u \in {\X}Q^{2n-1}.
\end{array} $$
\end{prop}
We note that there are three types of $\alpha$-Sasakian space
forms with respect to $sign(c+3{\alpha}^2)$ \cite {T}:
$$\begin{array}{l}
Type \, I: c+3{\alpha}^2 > 0;\\
[2mm]
Type \, II: c+3{\alpha}^2 = 0;\\
[2mm] Type \, III: c+3{\alpha}^2 < 0.
\end{array}$$

\section{K\"ahler-Lorentz manifolds with $B_0$-distributions}
\vskip 2mm Let $(M, h', J) \, (\dim M = 2n)$ be a complex manifold
with complex structure $J$ and indefinite Hermitian metric $h'$ of
signature $(2(n - 1), 2)$ and $\nabla'$ be the Levi-Civita
connection of $h'$. If $\nabla' J = 0,$ then $(M, h', J)$ is said
to be a {\it K\"ahler-Lorentz manifold}.

We consider K\"ahler-Lorentz manifolds $(M, h', J)$ with a
space-like $J$-invariant distribution $D$ of $\dim D = 2(n - 1).$
Then the orthogonal $J$-invariant two-dimensional distribution
$D^{\perp}$ is time-like.

Since our considerations are local, we can assume the existence of
a time-like unit vector field $\xi'$ on $M$ such that
$D^{\perp}(p) = span \{\xi', J\xi'\}$ at any point $p \in M $. We
denote by $\eta'$ and $\tilde \eta'$ the unit 1-forms
corresponding to $\xi'$ and $J\xi'$, respectively, i.e.
$$\eta'(X) = h'(\xi', X), \quad \tilde \eta' (X) = h'(J\xi', X)
= -\eta'(JX), \quad X \in {\X}M;$$
$$\| \eta' \| ^ 2 = \| \tilde \eta' \| ^ 2 =
\eta'(\xi') = \tilde\eta'(J\xi') = -1.$$

Then the space-like distribution $D$ is determined by the
conditions
$$D(p) = \{X \in T_pM \, \vert \, \eta'(X) = \tilde \eta'(X)
= 0 \}, \quad p \in M.$$

The Riemannian curvature tensor $R'$ of $\nabla'$ is determined as
in the previous section. We note that the Ricci tensor $\rho'$ and
the scalar curvature $\tau'$ of the metric $h'$ are given by
$$\rho'(Y,Z)=\sum_{i=1}^{2n}h'(e_i,e_i)R'(e_i,Y,Z,e_i),
\quad Y,Z \in {\X}M;$$
$$\tau'=\sum_{i=1}^{2n}h'(e_i,e_i)\rho'(e_i,e_i),$$
where $\{e_i\}, \, i=1,...,2n$ is an orthonormal basis for
$T_pM,\,p\in M$.

We also note that the tensor $h'^{\perp}=-(\eta'\otimes\eta'
+\tilde\eta'\otimes\tilde\eta')$ does not depend on the basis
$\{\xi',J\xi'\}$ of $D^{\perp}$. This tensor is negative definite
and it is the restriction of the metric $h'$ onto the distribution
$D^{\perp}$.

The K\"ahler form $\Theta$ of the structure $(h',J)$ is given by
$\Theta(X,Y)=h'(JX,Y), \; X,Y \in {\X}M$.

All directions in $D^{\perp}= span \{\xi',J\xi'\}$ have one and
the same Ricci curvature, which is denoted by $\sigma'$, i.e.
$$\sigma' =-\rho'(\xi',\xi')=-\rho'(J\xi',J\xi').\leqno{(3.1)}$$
The Riemannian sectional curvature of the distribution $D^{\perp}$
is denoted by $\varkappa'$, i.e.
$$\varkappa' =R'(\xi',J\xi',J\xi',\xi').\leqno{(3.2)}$$
Thus the structure $(h',J,D)$ gives rise to the functions
$\varkappa', \sigma'$ and $\tau'$.

Any vector field $X \in {\X}M$ is decomposable in a unique way as
follows:
$$X = x_0 - \tilde \eta'(X)J\xi' - \eta'(X)\xi',$$
where $x_0$ is the projection of $X$ into ${\X}D$.

As a rule, we use the following denotations for vector fields
(vectors):
$$X,Y,Z \in {\X}M \; (T_pM); \quad x_0,y_0,z_0 \in {\X}D \; (D(p)).$$

If $D^{\perp} = span \{\xi',J\xi'\}$, then the relative
divergences $div_0\xi'$ and $div_0J\xi'$ (the relative
codifferentials $\delta_0\eta'$ and $\delta_0\tilde\eta'$) of the
vector fields $\xi'$ and $J\xi'$ (of the 1-forms $\eta'$ and
$\tilde\eta'$) with respect to the space-like distribution $D$ are
introduced as in the definite case:
$$div_0\xi' = -\delta_0\eta' = \sum_{i = 1}^{2(n-1)}
(\nabla'_{e_i}\eta')e_i, \quad div_0J\xi' = - \delta_0\tilde\eta'
= \sum_{i = 1}^{2(n-1)}(\nabla'_{e_i}\tilde\eta')e_i,$$ where
$\{e_1,...,e_{2(n-1)}\}$ is an orthonormal basis of $D(p), \, p
\in M$.

The restriction of the metric $h'$ onto the distribution $\Delta$
$$\Delta(p) := \{X \in T_p M\; \vert\; \eta'(X) = 0\},\; p \in M,$$
perpendicular to $\xi'$, is of signature $(2(n-1),1)$.

The notion of a space-like $B_0$-distribution in a
K\"ahler-Lorentz manifold is introduced similarly to the definite
case:

\begin{defn}
Let $(M,h',J,D) \, (\dim \, M = 2n \geq 6)$ be a K\"ahler-Lorentz
manifold with $J$-invariant space-like distribution $D$
$(D^{\perp} = span\{\xi',J\xi'\})$. The distribution $D$ is said
to be a {\it $B_0$-distribution} if:
$$\begin{array}{l}
i) \displaystyle{ \quad \nabla' _{x_0} \xi' =
-\frac{k'}{2}\,x_0, \; k'\neq 0,}\quad x_0 \in D;\\
[2mm]
ii) \quad \nabla'_{J\xi'}\xi'={p^*}'J\xi';\\
[2mm]
iii)\quad \nabla' _{\xi'} \xi' = 0, \\
[2mm]
\end{array}\leqno{(3.3)}$$
where $k'$ and ${p^*}'$ are functions on $M$.
\end{defn}
Next we prove some properties of K\"ahler-Lorentz manifolds with
$B_0$-distribution.

\begin{lem}\label{L:3.1}
Let $(M,h',J,D) \, (\dim \, M = 2n \geq 6)$ be a K\"ahler-Lorentz
manifold with $B_0$-distribution $D$ $(D^{\perp} =
span\{\xi',J\xi'\})$. Then
$$dk'=-\xi'(k')\eta', \quad {p^*}'=\frac{\xi'(k')-k'^2}{k'}.$$
\end{lem}
{\it Proof}. The conditions (3.3) imply
$$\nabla'_X\xi'=-\frac{1}{2}k'\{X+\tilde\eta'(X)J\xi'
+\eta'(X)\xi'\}-{p^*}'\tilde\eta'(X)J\xi'. \leqno{(3.4)}$$ By
using (3.4) we find $d\tilde\eta'$ and after an exterior
differentiation we obtain the assertion of the lemma. \hfill{\bf
QED}

\begin{lem}\label{L:3.2}
Let $(M,h',J,D) \, (\dim \, M = 2n \geq 6)$ be a K\"ahler-Lorentz
manifold with $B_0$-distribution $D$ $(D^{\perp} =
span\{\xi',J\xi'\})$. Then
$$\begin{array}{ll}
\vspace{2mm} R'(X,Y)\xi'=&\displaystyle{\frac{1}{2}\left(
\xi'(k')-\frac{1}{2}k'^2
\right)\{\eta'(X)Y-\eta'(Y)X}\\
\vspace{2mm}
&+2h'(JX,Y)J\xi'-\tilde\eta'(X)JY+\tilde\eta'(Y)JX\}\\
\vspace{2mm}
&-\displaystyle{\frac{1}{k'}\,\xi'\left(\xi'(k')-\frac{1}{2}k'^2\right)
(\eta'\wedge\tilde\eta')(X,Y)J\xi';}\end{array}\leqno{(3.5)}$$
$$\varkappa'=-\frac{1}{k'}\,\xi'\left(\xi'(k')-\frac{1}{2}k'^2\right)
-2\left(\xi'(k')-\frac{1}{2}k'^2\right);\leqno{(3.6)}$$
$$\sigma'=-\frac{1}{k'}\,\xi'\left(\xi'(k')-\frac{1}{2}k'^2\right)-
(n+1)\left(\xi'(k')-\frac{1}{2}k'^2\right).\leqno{(3.7)}$$
\end{lem}
{\it Proof}. By using (3.4), we find immediately (3.5) and (3.6).
Taking a trace in (3.5), we have
$$\rho'(Y,\xi')=-\left[\frac{1}{k'}\,\xi'\left(\xi'(k')
-\frac{1}{2}k'^2\right)+(n+1)\left(\xi'(k')-\frac{1}{2}k'^2
\right)\right]\eta'(Y),\leqno{(3.8)}$$ which implies (3.7).
\hfill{\bf QED} \vskip 2mm The equality (3.8) shows that every
unit vector in $D^{\perp}(p)$ is an eigen vector of the Ricci
operator $\rho'$ with one and the same eigen value $\sigma'(p)$.

If $x_0$ is a unit vector in $D(p)$, then the Riemannian sectional
curvature of $span\{x_0,\xi'\}$ may only depend on the point $p
\in M$:
$$-R'(x_0,\xi',\xi',x_0)=\frac{\sigma'-\varkappa'}{2(n-1)}.
\leqno{(3.9)}$$

The first step in the study of K\"ahler-Lorentz manifolds with
$B_0$-distributions is to describe the flat case.

Let $(M,h',J,D)\; (\dim M \geq 6)$ be a flat K\"ahler-Lorentz
manifold with $B_0$-distribution $D$ $(D^{\perp} =
span\{\xi',J\xi'\})$. Then Lemma \ref{L:3.2} implies that
$$\xi'(k')=\frac{1}{2}\,k'^2.\leqno {(3.10)}$$
Taking into account Lemma \ref{L:3.1} it follows that
$${p^*}'=-\frac{1}{2}\,k'.\leqno {(3.11)}$$
Then (3.4) in view of (3.10) and (3.11) implies that
$$ \nabla'_x\xi' = -\frac{1}{2}\,k'\,x,\quad h'(x,\xi')=0.\leqno {(3.12)}$$
Hence the integral submanifolds $Q^{2(n-1)}_1$ of the distribution
$\Delta$, perpendicular to $\xi'$, are totally umbilic
submanifolds of $M$ with time-like normals $\xi'$.

Let $(\mathbb{C}^n = \{\textbf{Z}=(z^1,..., z^{n-1}; z^n)\}, J)$
be the standard $n$-dimensional complex vector space with complex
structure $J$ and $h'$ be the K\"ahler metric of signature
$(2(n-1),2)$, defined by
$$h'(\textbf{Z}, \textbf{Z})= \vert z^1 \vert^2 +...
+ \vert z^{n-1} \vert^2 - \vert z^n \vert^2.$$ We call $h'$ the
{\it canonical flat K\"ahler-Lorentz metric} and $(\mathbb{C}^n,
h', J)=(\mathbb{R}^{2(n-1)}_2, h', J)$ the {\it canonical flat
K\"ahler-Lorentz manifold}.

Next we describe the $B_0$-distributions in $(\mathbb{C}^n, h',
J)$.

Let $D$\,$(D^{\perp} = span \{\xi', J\xi'\})$ be a
$B_0$-distribution in $(\mathbb{C}^n, \,h', J).$ According to
Definition 3.1 $\xi'$ is a time-like geodesic vector field with
respect to the flat Levi-Civita connection $\nabla'$ of $h'$. Then
the integral curves of $\xi'$ are straight lines. Since $h'$ is
flat, then the integral submanifolds $Q^{2(n-1)}_1$ of the
distribution $\Delta$, perpendicular to $\xi'$, are totally
umbilical with time-like normals $\xi'$. Applying the standard
theorem for totally umbilical submanifolds (with time-like
normals) of the manifold $(\mathbb{C}^n, h', J)$, we obtain that
$Q^{2(n-1)}_{1}$ is locally a part of a hypersphere
$H_{1}^{2(n-1)}(\textbf{Z}_0, r): h'(\textbf{Z}-\textbf{Z}_0,
\textbf{Z}-\textbf{Z}_0) = - r^2,\quad r>0.$ All these
hyperspheres are orthogonal to the integral curves of $\xi'$, i.e.
$Q^{2(n-1)}_1$ are the concentric hyperspheres
$$H_{1}^{2(n-1)}(\textbf{Z}_0, r):
h'(\textbf{Z}-\textbf{Z}_0, \textbf{Z}-\textbf{Z}_0) = - r^2, \,
\textbf{Z}_0 = const.$$

Choosing $\textbf{Z}_0$ at the origin $O$ of ${\C}^n$, we obtain
\vskip 2mm

{\it Canonical example of a flat K\"ahler-Lorentz manifold with
$B_0$-distribution:}
$$({\T}^{n-1}_{1}, h', J, D),$$
where ${\T}^{n-1}_1$ is the time-like domain in $\mathbb{C}^n$
$$ {\T}^{n-1}_{1} = \{ \textbf{Z}\in {\C}^{n} \;
\vert \; h'(\textbf{Z}, \textbf{Z}) < 0 \}$$ and
$$\xi' = \frac{\textbf{Z}}{\sqrt{-h'(\textbf{Z},\textbf{Z})}},
\quad \textbf{Z} \in \mathbb{T}^{n-1}_1.$$

Now let $(M, h', J, D)\, (\dim M = 2n \geq 6)$ be a flat
K\"ahler-Lorentz manifold with $B_0$-distribution $D \,(D^{\perp}
= span \{\xi',J\xi'\}).$ Since the Levi-Civita connection
$\nabla'$ of $h'$ is flat and $\nabla' J = 0,$ then there exists a
local holomorphic isometry $\phi$ of $(M, h', J)$ onto
$({\C}^n,h',J)$. Since $\phi$ transforms the $B_0$-distribution
$D$ into a $B_0$-distribution, then we have

\begin{prop} \label{P:3.1}
Any flat K\"ahler-Lorentz manifold with $B_0$-distribution is
locally equivalent to the canonical example
$(\mathbb{T}^{n-1}_1,h',J,D)$.
\end{prop}

In order to make computations in local holomorphic coordinates we
need some formulas concerning the structures on
$\mathbb{T}^{n-1}_{1}$.

Let in $\mathbb{C}^n=\{\textbf{Z}=(z^1,...z^n)\} \; (n \geqq 2)$
\; $\partial_{\alpha} := \displaystyle {\frac{\partial}{\partial
z^{\alpha}},\;
\partial_{\bar\alpha} := \frac{\partial}{\partial z^{\bar\alpha}}
= \overline{\frac{\partial}{\partial z^{\alpha}}}},\; \alpha =
1,..., n.$ Further the indices $\alpha, \beta,...$ will run over
$1,..., n$.

The canonical flat K\"ahler-Lorentz metric $h'$ has the following
local components:
$$h'_{\alpha\bar\beta}=\left \{
\begin{array}{ll}
\vspace{2mm}
\displaystyle{\; \; \; \frac{1}{2}}&\alpha=\beta=1,...,n-1;\\
\vspace{2mm}
\displaystyle{-\frac{1}{2}}&\alpha=\beta=n;\\
\vspace{2mm} \; \; \;0&\alpha\neq\beta.
\end{array}\right. $$
Then
$$h'(\textbf{Z},\textbf{Z})=\vert z^1 \vert^2+...+\vert z^{n-1} \vert ^2-
\vert z^n \vert ^2 = 2
h'_{\alpha\bar\beta}\,z^{\alpha}z^{\bar\beta},$$ where the
summation convention is assumed.

The distance function $-r^2=h'(\textbf{Z},\textbf{Z})$ in the
domain $\mathbb{T}^{n-1}_{1}$ is given by
$$-r^2=2h'_{\alpha\bar\beta}z^{\alpha}z^{\bar\beta}<0,\quad r>0.
\leqno{(3.13)}$$ The vector field
$\xi'=\displaystyle{\frac{1}{r}\,\textbf{Z}}$ at the point $p \in
{\T}^{n-1}_{1}$ with position vector $\textbf{Z}$ has local
components
$$\eta'^{\alpha}=\frac{1}{r}\,\delta^{\alpha}_{\sigma}z^{\sigma},$$
where $\delta_{\alpha}^{\sigma}$ are the Kronecker's deltas.

Taking into account (3.13), we find the local components of the
corresponding 1-form $\eta'$:
$$\eta'_{\alpha}=\eta'^{\bar\sigma}h'_{\alpha\bar\sigma}=
\frac{1}{r}\,h'_{\alpha\bar\beta}z^{\bar\beta}=-r_{\alpha}.
\leqno{(3.14)}$$ Hence
$$\begin{array}{l}
\vspace{2mm}
\displaystyle{\eta'=-dr, \quad \xi' = \frac{d}{dr}};\\
\vspace{2mm} \displaystyle{\eta'(\xi')=h'(\xi',\xi')=
\frac{1}{r^2}\,h'(\textbf{Z},\textbf{Z})=-1.}
\end{array}
\leqno{(3.15)}$$

By differentiating (3.13) we obtain
$$h'_{\alpha\bar\beta}=\frac{1}{2}\,\partial_{\alpha}
\partial_{\bar\beta}(-r^2).$$

On the other hand, differentiating (3.14), we have
$$\partial_{\bar\beta}\eta'_{\alpha}=\nabla'_{\bar\beta}\eta'_{\alpha}=
\frac{1}{r}\,(h'_{\alpha\bar\beta}+\eta'_{\alpha}\eta'_{\bar\beta}).$$

\begin{prop}\label{P:3.2}
Let $f(t),\, t<0$ be a real $\mathcal{C}^{\infty}$-function
satisfying the inequalities:
$$f'(t)>0, \quad f'(t)+tf''(t)<0.$$

Then
$$g_{\alpha\bar\beta}=\partial_{\alpha}\partial_{\bar\beta} f(-r^2)$$
are the local components of a K\"ahler metric g.
\end{prop}
{\it Proof}. By using (3.13) and (3.15), we calculate
$$\partial _{\bar\beta} f(-r^2)=2f'h'_{\alpha\bar\beta}z^{\alpha}.$$
Differentiating the last equality, we find
$$g_{\alpha\bar\beta}=\partial_{\alpha}\partial_{\bar\beta}f(-r^2)=
2f'h'_{\alpha\bar\beta}+4r^2f''\eta'_{\alpha}\eta'_{\bar\beta}.$$
Hence,
$$g=2f'h'+2r^2f''(\eta'\otimes\eta'+\tilde\eta'\otimes\tilde\eta').
\leqno{(3.16)}$$

Now, let $p\in {\T}_1^{n-1}$ and $T_p({\T}_1^{n-1})= (D(p)\oplus
D^{\perp}(p))$. The equality (3.16) implies that
$$g(x_0,x_0)=2f'h'(x_0,x_0),\quad x_0\in D(p);\leqno{(3.17)}$$
$$g(\xi',\xi')=g(J\xi',J\xi')=-2(f'+(-r^2)f'').\leqno{(3.18)}$$
The first condition of the proposition and (3.17) imply that the
restriction of $g$ onto $D$ is positive definite. The second
condition of the proposition and (3.18) give that the restriction
of $g$ onto $D^{\perp}$ is also positive definite. Hence $g$ is a
positive definite metric. Since $g=\partial\bar\partial f(-r^2)$,
then $g$ is a K\"ahler metric. \hfill{\bf QED}

\section{K\"ahler manifolds of quasi-constant holomorphic
sectional curvatures with $a+k^2 < 0$}

In this section we prove the main theorem, which clarifies the
connection between the K\"ahler metrics introduced in Section 3
and a class of K\"ahler manifolds of quasi-constant holomorphic
sectional curvatures.

Let $(M, g, J, D)\; (\dim M = 2n \geq 6)$ be a K\"ahler manifold
with $B_0$-distribution $D\;(D^{\perp} = span \{\xi, J\xi\})$ with
functions $k, p^*$, given by (2.2).

If $u, v$ are proper $\mathcal{C}^{\infty}$-functions of the
distribution $\Delta$ (cf \cite {GM2}), i.e. $du = \xi(u)\,\eta,\;
dv = \xi(v)\,\eta,$ we consider the metric
$$h' = e^{2u}\left(g - (e^{2v} + 1)(\eta\otimes\eta +
\tilde\eta\otimes\tilde\eta) \right),\leqno (4.1)$$ which is
positive definite on $D$ and negative definite on $D^{\perp}$.

\begin{lem}\label{L:4.1}
Let $(M, g, J, D)\; (\dim M = 2n \geq 6)$ be a K\"ahler manifold
with $B_0$-distribution $D \;(D^{\perp} = span \{\xi, J\xi\})$.
Then the metric $h'$, given by $(4.1)$ is K\"ahler-Lorentz if and
only if
$$\displaystyle{\xi(u) = -\frac{k(e^{2v} + 1)}{2}.}\leqno (4.2)$$
\end{lem}

{\it Proof.} From (4.1) we find the K\"ahler form $\Theta$ of the
metric $h'$:
$$\Theta = e^{2u}\left(\Omega - (e^{2v} + 1)\eta\wedge\tilde\eta
\right).$$ The last equality, (2.1) and (2.2) imply that
$$d\Theta = e^{2u}\left(2\xi(u) + k(e^{2v} + 1)\right)
\eta\wedge\Omega,$$ which implies the assertion of the lemma.
\hfill{\bf QED} \vskip 2mm We set
$$\xi' = e^{-(u+v)}\,\xi, \quad \eta' = -e^{u+v}\,\eta.\leqno  (4.3)$$
Then $\eta'$ is the 1-form corresponding to $\xi'$ with respect to
$h'$ and $\eta'(\xi')=-1$.

\begin{lem}\label{L:4.2}
Let $(M, g, J, D)\; (\dim M = 2n \geq 6)$ be a K\"ahler manifold
with $B_0$-distribution $D \;(D^{\perp} = span \{\xi, J\xi\}).$ If
$$h' = e^{2u}\left(g - (e^{2v} + 1)(\eta\otimes\eta +
\tilde\eta\otimes\tilde\eta) \right),$$ where
$$dv = \xi(v)\,\eta, \quad du = -\frac{k(e^{2v} + 1)}{2}\,\eta$$
and
$$\xi' = e^{-(u+v)}\,\xi, \quad \eta' = -e^{u+v}\,\eta,$$
then $(M, h', J, D)\; (D^{\perp} = span \{\xi', J\xi'\})$ is a
K\"ahler-Lorentz manifold with space-like $B_0$-distribution $D$.
\end{lem}

{\it Proof.} Let $\nabla', \nabla$ be the Levi-Civita connections
of the metrics $h', g,$ respectively. Then
$$\begin{array}{ll}
\vspace{2mm} \nabla'_XY=&\nabla_XY + \xi(u)\{\eta(X)Y + \eta(Y)X +
\tilde\eta(X)JY + \tilde\eta(Y)JX \} \\
\vspace{2mm}
& + \xi(v-u)\{[\eta(X)\eta(Y) - \tilde\eta(X)\tilde\eta(Y)]\xi \\
\vspace{2mm} & +
[\eta(X)\tilde\eta(Y)+\tilde\eta(X)\eta(Y)]J\xi\}, \quad X,Y \in
{\X}M. \end{array}\leqno{(4.4)}$$ From (4.4) it follows that
$$\begin{array}{ll}
\vspace{2mm} \nabla'_X\xi'=&\displaystyle{e^{-(u+v)}(\xi(u) +
\frac{k}{2})
\left[X-\eta(X)\xi-\tilde\eta(X)J\xi\right]}\\
\vspace{2mm} & +
e^{-(u+v)}\left(\xi(u+v)-p^*\right)\tilde\eta(X)J\xi, \quad X\in
{\X}M.\end{array} $$

The above equality can be written in the form
$$\nabla'_X\xi'=-\frac{k'}{2}[X+\eta'(X)\xi'+\tilde\eta'(X)J\xi']
-{p^*}'\tilde\eta'(X)J\xi', \quad X\in {\X}M,\leqno (4.5)$$ where
$$k'=-2e^{-(u+v)}(\xi(u)+\frac{k}{2}),\leqno (4.6)$$
$${p^*}' = e^{-(u+v)}(\xi(u+v) - p^*),$$
i.e. $D$ is a space-like $B_0$-distribution with functions $k'$
and ${p^*}'$.\hfill{\bf QED} \vskip 2mm

Because of (4.2)
$$\xi(u) + \frac{k}{2} = -\frac{1}{2}e^{2v}k.$$
Then (4.6) gives the following relation between $k'$ and $k$:
$$ k' = e^{v-u}\,k.\leqno (4.7)$$
\vskip 2mm Let the tensors $\pi', \Phi'_1, \Phi'_2, \Phi' =
\Phi'_1 + \Phi'_2$ and $\Psi'$ of type (1,3) with respect to the
structure $(h', \xi', \eta')$ be determined as in (2.3), (2.4) and
(2.5). If $g$ and $h'$ are related as in Lemma \ref{L:4.2}, then
$$\begin{array}{l}
\vspace{2mm}
\pi'+ 2\Phi' + \Psi'=e^{2u}(\pi - 2\Phi + \Psi),\\
\vspace{2mm}
\Phi'_1+\frac{1}{2}\Psi'=-e^{2u}(\Phi_1-\frac{1}{2}\Psi),\\
\vspace{2mm}
\Phi'_2+\frac{1}{2}\Psi'=e^{2(u+v)}(\Phi_2-\frac{1}{2}\Psi),\\
\vspace{2mm} \Psi'=-e^{2(u+v)}\Psi.\end{array}\leqno (4.8)$$

\begin{prop}\label{P:4.1}
Let $(M, g, J, D)\;(\dim M = 2n \geq 6)$ be a K\"ahler manifold of
quasi-constant holomorphic sectional curvatures with
$B_0$-distribution $D \;(D^{\perp} = span \{\xi, J\xi\})$ and
$$a + k^2 < 0.$$
If the structure $(h', \xi', \eta')$ is determined as in Lemma
$\ref{L:4.2}$ by the proper function
$$e^{2v} = -\frac{a+k^2}{k^2},$$
then $h'$ is a flat K\"ahler-Lorentz metric.
\end{prop}
{\it Proof.} By direct computations from (4.4) in view of (2.3),
(2.4) and (2.5) we find
$$\begin{array}{ll}
R'-R=&-2k\xi(u)(\pi - 2\Phi +\Psi)\\
[2mm] &\displaystyle{-4k\xi(v)(\Phi_1-\frac{1}{2}\Psi)
-4(\xi^2(u)-p^*\xi(u))(\Phi_2-\frac{1}{2}\Psi)}\\
[3mm] &-(\xi^2(u+v)-p^*\xi(u+v))\Psi.\end{array}\leqno{(4.9)}$$
Taking into account that $R = a\pi + b\Phi + c\Psi$ and (4.8), we
obtain from (4.9) the curvature tensor $R'$ of $h'$ in the form
$$R' = A(\pi'+2\Phi'+\Psi') + B_1(\Phi_1'+\frac{1}{2}\Psi')
+B_2(\Phi_2'+\frac{1}{2}\Psi') + C\Psi',\leqno(4.10)$$ where
$$\begin{array}{ll}
&e^{2u}A=a-2k\xi(u), \quad e^{2(u+v)}C=-(a+b+c)+{\xi}^2(u+v) -p^*\xi(u+v),\\
[2mm] &e^{2u}B_1=-(2a+b)+4k\xi(v),\quad
e^{2(u+v)}B_2=2a+b-4(\xi^2(u)-p^*\xi(u)).
\end{array}\leqno(4.11)$$
Taking into account (4.7), (4.10) and (4.2), we find
$$e^{2u}(A-{k'}^2)=a+k^2.\leqno (4.12)$$
Then (4.12) and (4.7) imply
$$e^{2u}A=e^{2v}k^2+a+k^2.$$
Under the conditions of the proposition we obtain $A=0$ and
$\xi(u) = \displaystyle{\frac{a}{2k}}.$

Differentiating the equality $\displaystyle{e^{2v}=
-\frac{a+k^2}{k^2}}$, because of (2.6), we obtain
$$\xi(k) + \frac{1}{2}k^2 + k\xi(v) = 0.\leqno (4.13)$$
On the other hand, $\xi' = e^{-(u+v)}\,\xi$ and (4.7) imply
$$\xi(k) + \frac{1}{2}k^2 + k\xi(v) = e^{2u}(\xi'(k') -
\frac{1}{2}{k'}^2).$$ Thus, from the equality (4.13), we get
$$\xi'(k') = \frac{1}{2}\,{k'}^2.$$
Now from (3.6) and (3.7) it follows that $\varkappa' = \sigma' =
0.$

Replacing into (4.10) the quadruples \;$\xi, x_0, x_0, \xi$
\;and\; $x_0, \xi, \xi, x_0,$\;where \;$h'(x_0, x_0)=1,$\; in view
of (3.9), we obtain
$$0=\frac{\sigma'-\varkappa'}{2(n-1)}=\frac{1}{8}B_1=\frac{1}{8}B_2.$$
Replacing into (4.10) the quadruple\; $\xi, J\xi, J\xi, \xi,$\; we
get
$$0=\varkappa'=C,$$
i.e. $R'=0.$ \hfill{\bf QED} \vskip 2mm

Let now $(M, h', J, D)\; (\dim M = 2n \geq 6)$ be a
K\"ahler-Lorentz manifold with space-like $B_0$-distribution
$D\;(D^{\perp} = span \{\xi', J\xi'\})$ with functions $k'$ and
${p^*}'$, determined in Lemma \ref{L:3.1}.

If $u,v$ are proper ${\mathcal{C}}^{\infty}$-functions of the
distribution $\Delta$, i.e. $du=-\xi'(u)\eta', \;
dv=-\xi'(v)\eta',$ we consider the metric
$$g=e^{-2u}(h'+(e^{-2v}+1)(\eta'\otimes\eta'+
\tilde\eta'\otimes\tilde\eta')). \leqno (4.14)$$ Taking into
account (3.4), analogously to Lemma \ref{L:4.1}, we have

\begin{lem}\label{L:4.3}
Let $(M, h', J, D)\; (\dim M = 2n \geq 6)$ be a K\"ahler-Lorentz
manifold with space-like $B_0$-distribution $D \;(D^{\perp} = span
\{\xi', J\xi'\})$. Then the metric $g$, given by $(4.14)$ is
K\"ahler if and only if
$$\displaystyle{\xi'(u) =-\frac{k'(e^{-2v} + 1)}{2}.}$$
\end{lem}

Further we set $\xi=e^{u+v}\xi',\; \eta=-e^{-(u+v)}\eta'.$
Analogously to (4.7) we have
$$k=e^{u-v}k'.\leqno (4.15)$$

\begin{lem}\label{L:4.4}
Let $(M, h', J, D)\; (\dim M = 2n \geq 6)$ be a K\"ahler-Lorentz
manifold with space-like $B_0$-distribution $D \;(D^{\perp} = span
\{\xi', J\xi'\}).$ If
$$g = e^{-2u}\left(h' + (e^{-2v} + 1)(\eta'\otimes\eta' +
\tilde\eta'\otimes\tilde\eta') \right),$$ where
$$dv = -\xi'(v)\,\eta', \quad du = \frac{k'(e^{-2v} + 1)}{2}\,\eta'$$
and
$$\xi = e^{u+v}\,\xi', \quad \eta = -e^{-(u+v)}\,\eta',$$
then $(M, g, J, D)\; (D^{\perp} = span \{\xi, J\xi\})$ is a
K\"ahler manifold with $B_0$-distribution $D$.
\end{lem}
\vskip 2mm
\begin{prop}\label{P:4.2}
Let $({\mathbb{T}}^{n-1}_1, h', J, D)\;(\dim {\mathbb{T}}^{n-1}_1
= 2n \geq 6)$ be the canonical example of a flat K\"ahler-Lorentz
manifold. If the structure $(g, \xi, \eta)$ is determined as in
Lemma $\ref{L:4.4}$, then $g$ is a K\"ahler metric of
quasi-constant holomorphic sectional curvatures and $a+k^2 < 0$.
\end{prop}

{\it Proof.} Taking into account (4.14) we find the relation (4.4)
between the Levi-Civita connections $\nabla'$ and $\nabla$ of $h'$
and $g$, respectively. Then the corresponding relation between the
curvature tensors $R'$ and $R$ is given by (4.9). Since $R'=0$,
then (4.9) gives the tensor $R$ in the form
$$R=A^*(\pi-2\Phi+\Psi)+B^*_1(\Phi_1-\frac{1}{2}\Psi)
+B^*_2(\Phi_2-\frac{1}{2}\Psi)+C^*\Psi.$$ Replacing the quadruples
\;$\xi, x_0, x_0, \xi;\; x_0, \xi, \xi, x_0,$\;where \;$g(x_0,
x_0)=1$,\; in the last equality, we get
$$\frac{1}{8}B^*_1=R(\xi, x_0, x_0, \xi)=R(x_0, \xi, \xi, x_0)
=\frac{1}{8}B^*_2.$$ Hence the curvature tensor $R$ has the form
$R=a\pi+b\Phi+c\Psi$, i.e. the metric $g$ is of quasi-constant
holomorphic sectional curvatures.

To prove $a+k^2 < 0$, we consider (4.9). Since $a=2k\xi(u),\;
\xi=e^{u+v}\xi',$\; in view of (4.15) and Lemma \ref{L:4.3}, we
find
$$a=-e^{2u}{k'}^2(e^{-2v}+1)=-k^2-e^{2u}{k'}^2.$$
Hence $a+k^2 < 0$. \hfill{\bf QED} \vskip 2mm
\begin{thm}\label{T:4.1}
Any K\"ahler metric $g=\partial\bar\partial f(-r^2),$ $-r^2$ being
the time-like distance function in ${\T}^{n-1}_1$ $(n \geq 3),$ is
of quasi-constant holomorphic sectional curvatures and function
$a+k^2 < 0.$

Conversely, every K\"ahler manifold $(M,g,J,D)\,(\dim M=2n \geq
6)$ of quasi-constant holomorphic sectional curvatures with
$B_0$-distribution satisfying the condition $a+k^2 < 0$ is locally
equivalent to $({\T}^{n-1}_1,g,J,D)$ with the canonical
$B_0$-distribution and $g=\partial\bar\partial f(-r^2)$.
\end{thm}

{\it Proof.} Let the K\"ahler metric $g$ be given as in (3.16).
Putting
$$e^{-2u}=2f',\quad  e^{-2v}+1=\frac{r^2f''}{f'},$$
we calculate
$$\xi'(u) =\frac{du}{dr}=\frac{rf''}{f'}=\frac{1}{r}(e^{-2v}+1)
=-\frac{k'(e^{-2v} + 1)}{2}.$$ Then we can apply Proposition \ref
{P:4.2} and conclude that the structure $(g, J, \xi, \eta)$ is of
quasi-constant holomorphic sectional curvatures with
$B_0$-distribution and function \;$a+k^2 < 0.$

For the inverse, let $(M, g, J, D)\; (\dim M = 2n \geq 6)$ be a
K\"ahler manifold with $B_0$-distribution and function $a+k^2 <
0.$ We construct the metric $h'$ as in Lemma \ref{L:4.2} by the
proper function $\displaystyle{e^{2v}=-\frac{a+k^2}{k^2}}$.
Applying Proposition \ref {P:4.2} we obtain the K\"ahler metric
$h'$ is flat and the given manifold is locally equivalent to the
canonical flat K\"ahler-Lorentz manifold $({\T}_1^{n-1},h',J,D)$.

Further we write the equality (4.1) in the form
$$g = e^{-2u}\left(h' + (e^{-2v} + 1)(\eta'\otimes\eta' +
\tilde\eta'\otimes\tilde\eta') \right)\leqno(4.16)$$ and put
$$f(-r^2)=\frac{1}{2}\int{e^{-2u}d(-r^2)}.\leqno{(4.17)}$$

From (4.17) we have $e^{-2u}=2f'$. Using Lemma \ref{L:4.4} we find
$\displaystyle{\xi'(u)=-\frac{k'(e^{-2v}+1)}{2}}$\; and\;
$\displaystyle{\xi'=\frac{d}{dr},}$ \;
$\displaystyle{k'=-\frac{2}{r}}.$\; Then
$\displaystyle{e^{-2v}+1=\frac{r^2f''}{f'}}$ and (4.16) becomes
$$g=2f'(h'+\frac{r^2f''}{f'}(\eta'\otimes\eta' +
\tilde\eta'\otimes\tilde\eta')).$$ Hence $g=\partial\bar\partial
f(-r^2)$ with potential function (4.17). \hfill{\bf QED} \vskip
2mm

As an application of Theorem \ref{T:4.1} we shall find the
K\"ahler metrics of constant holomorphic sectional curvatures,
defined in the manifold $(\mathbb{T}^{n-1}_1, h', J, D)$ by the
condition $g=\partial\bar\partial f(-r^2).$

Let $g$ be a metric given by (4.16). Then (4.9) gives the relation
between the curvature tensor $R$ of $g$ and the tensor $R'=0$ of
$h'$ in $\mathbb{T}^{n-1}_1.$ Since the coefficients $A, B_1, B_2,
C$ in (4.10) are all zero, then (4.11) implies
$$\begin{array}{l}
\vspace{2mm}
a=2k\xi(u);\\
\vspace{2mm}
2a+b=4k\xi(v)=4[{\xi}^2(u)-p^*\xi(u)]\\
\vspace{2mm} a+b+c=\xi^2(u+v)-p^*\xi(u+v).
\end{array} \leqno (4.18)$$

Since $D$ is a $B_0$-distribution, then $g$ is a K\"ahler metric
of constant holomorphic sectional curvatures if and only if $b=0$.
Because of $\xi=e^{u+v}\xi'=e^{u+v} \displaystyle{\frac{d}{dr}}$
and (4.18) the condition $b=0$ is equivalent to the relation
$$\frac{du}{dr}=\frac{dv}{dr}.\leqno {(4.19)}$$

Further, taking into account (4.18) and (2.2), we obtain
successively
$$k\xi(u)=\xi^2(u)-p^*\xi(u)=\xi^2(u)+\frac {\xi(k)+k^2}{k^2}\xi(u),$$
which in view of (4.19) and the relation
$k=e^{u-v}k'=\displaystyle{-\frac{2e^{u-v}}{r}}$ implies that
$$\frac{d^2 u}{dr^2}+2\left(\frac{du}{dr}\right)^2 -
\frac{1}{r} \frac{du}{dr} =0. \leqno (4.20)$$

Solving (4.20), we find
$$e^{2u}=e^{2u_0} \vert r^2+a_0 \vert, \quad a_0=const,\quad u_0=const.
\leqno{(4.21)}$$

Since $a+k^2<0$, then $a<0$ and the equality
$a=2k\xi(u)=\displaystyle{-\frac{4}{r}e^{2u}\frac{du}{dr}=
-\frac{4e^{2u}}{r^2+a_0}}$ implies that $r^2+a_0>0$.

On the other hand, using the relation (4.2), we find
$e^{-2v}=\displaystyle{-\frac{a_0}{r^2+a_0}>0}$ and $a_0<0$.
Putting $a_0=-r_0^2$, we have
$$e^{-2v}=\frac{r_0^2}{r^2-r_0^2}.\leqno{(4.22)}$$

Finally, the equality $a=2k\xi(u)$ gives that $\displaystyle{
e^{-2u_0}=-\frac{4}{a}}$.

Now, from (4.21) and (4.22) we obtain \vskip 2mm {\it Examples of
K\"ahler space forms with $B_0$-distribution and $a + k^2 < 0$}:
\vskip 2mm All K\"ahler metrics $g$ of constant holomorphic
sectional curvatures $a<0,$ given in ${\T}^{n-1}_1$ by (4.16), are
$$g=-\frac{4}{a(r^2-r_0^2)}\left (h'+\frac{r^2}{r^2-r_0^2}\,
(\eta'\otimes\eta'+\tilde\eta'\otimes\tilde\eta')\right ), \,
r_0=const >0, \quad r>r_0.\leqno (4.23)$$

The potential function of the above metrics up to a constant is
$$f(-r^2)=\frac{2}{a} \ln (r^2-r_0^2), \, r_0=const >0, \quad r>r_0.$$

Hence
$$g=\frac{2}{a}\,\partial\bar\partial\,\ln (r^2-r_0^2),\, r_0>0,
\quad r>r_0.$$

One of these metrics is most remarkable:
$$g=\frac{4}{r^2-1}\left(h'+\frac{r^2}{r^2-1}
(\eta'\otimes\eta'+\tilde\eta'\otimes\tilde\eta')\right ),\quad
r>1. \leqno{(4.24)}$$ This metric is defined in the hyperbolic
unit "disc" ${\D}_1^{n-1}(1): h'(\textbf{Z},\textbf{Z})<-1$ and is
of constant holomorphic sectional curvatures $a=-1$.

\section{The geometric meaning of the function
$a+k^2$ in K\"ahler manifolds of quasi-constant holomorphic
sectional curvatures} \vskip 2mm Let $(M, g, J, D) \;(\dim M = 2n
\geq 6)$ be a K\"ahler manifold of quasi-constant holomorphic
sectional curvatures with $B_0$-distribution $D(p) \;(D^{\perp}(p)
= span \{\xi, J\xi\}),\, p \in M.$

In this section we study the geometric structure of the integral
submanifolds of the distribution
$$\Delta(p) = \{X \in T_pM \, \vert \, \eta(X) = 0\}, \, p\in M.$$

Because of (2.1), we have
$$ d\eta=0; \quad d\tilde\eta=k\Omega+\frac{1}{k}\eta\wedge\tilde\eta.
\leqno (5.1)$$

Let $Q^{2n-1}$ be an arbitrary integral submanifold of the
distribution $\Delta$ and $\xi$ be the unit vector field, normal
to $Q^{2n-1}.$ Applying the Weingarten and Gauss equations to the
submanifolds $Q^{2n-1}$, we have
$$\nabla_x\xi=\frac{k}{2}\,x+\frac{1}{k}\left(\xi(k)+\frac{k^2}{2}\right)
\tilde\eta(x)J\xi,\quad x \in {\X}\Delta;\leqno (5.2)$$
$$\nabla_xy=\mathcal{D}_xy+h(x, y)\xi, \quad x, y \in {\X}\Delta,
\leqno (5.3)$$ where $\mathcal{D}$ is the induced Levi-Civita
connection and $h$ is the second fundamental tensor on $Q^{2n-1}.$

According to (2.2) $k=const$ on $Q^{2n-1}.$ From (5.2) it follows
that
$$h=-\frac{k}{2}\;g-\frac{1}{k}\left(\xi(k)+\frac{k^2}{2}\right )
\tilde\eta\otimes\tilde\eta.$$

The standard almost contact Riemannian structures $(g, \varphi,
\tilde\xi, \tilde\eta)$ induced on the manifold $Q^{2n-1}$ are
\,\cite{T1},\,\cite {T2}:
$$\begin{array}{l}
\tilde\xi:=J\xi; \quad \tilde\eta=g(x,\tilde\xi),\\
[2mm] \varphi x:=Jx+\tilde\eta(x)\xi, \quad x \in {\X}\Delta.
\end{array}\leqno(5.4)$$

Taking into account (5.2), in view of (2.7), we find
$$\mathcal{D}_x\tilde\xi=\frac{k}{2}\;\varphi x, \quad x \in {\X}\Delta,
\leqno (5.5)$$
$$(\mathcal{D}_x \varphi)(y)=\frac{k}{2}\left (\tilde\eta(y)x-g(x, y)
\tilde\xi \right ),\quad x, y \in {\X}\Delta. \leqno (5.6)$$

According to (5.5) and (5.6), any integral submanifold $Q^{2n-1}$
of the distribution $\Delta$ is an $\alpha$-Sasakian manifold with
$\alpha=\displaystyle{\frac{k}{2}}$.

More precisely we have

\begin{prop}\label{P:5.2}
Let $(M, g, J, D)\; (\dim M = 2n \geq 6)$ be a K\"ahler manifold
of quasi-constant holomorphic sectional curvatures with
$B_0$-distribution $D \;(D^{\perp}=span\{\xi, J\xi\}).$

Then any integral submanifold $Q^{2n-1}$ of the distribution
$\Delta$
is a $\displaystyle{\frac{k}{2}}$\,-Sasakian space form\\
of type $\left \{\begin{array}{l} \vspace{1mm}
I,\\
\vspace{1mm}
II,\\
\vspace{1mm} III,
\end{array} \right .$
if and only if $\left \{
\begin{array}{l}
\vspace{1mm}
a+k^2>0,\\
\vspace{1mm}
a+k^2=0,\\
\vspace{1mm} a+k^2<0,
\end{array} \right .$ respectively.
\end{prop}
{\it Proof.} From (5.3), (5.2), (5.4) and (5.6) we find the
relation between the curvature tensors $R$ and $K$ of $M^{2n}$ and
$Q^{2n-1}$, respectively:
$$\begin{array}{ll}
R(x,y,z,u)=&K(x,y,z,u)-\displaystyle{\frac{1}{4}}k^2
[g(y,z)g(x,u)-g(x,z)g(y,u)]\\
[2mm] &-\displaystyle{\frac{1}{2}\left(\xi(k)+\frac{1}{2}k^2\right
)}
[g(y,z)\tilde\eta(x)\tilde\eta(u)+g(x,u)\tilde\eta(y)\tilde\eta(z)\\
[2mm]
&-g(x,z)\tilde\eta(y)\tilde\eta(u)-g(y,u)\tilde\eta(x)\tilde\eta(z)],
\quad x,y,z,u \in {\X}\Delta.
\end{array} \leqno (5.7)$$
Since $(M, g, J, D)$ is of quasi-constant holomorphic sectional
curvatures, then
$$R=a\pi+b\Phi+c\Psi. \leqno (5.8)$$
Taking into account (5.8), the equality (5.7) becomes
$$\begin{array}{ll}
K(x,y,z,u)=&\displaystyle{\frac{a+k^2}{4}}[g(y,z)g(x,u)-g(x,z)g(y,u)]\\
[2mm] &+\displaystyle{\frac{a}{4}}[g(\varphi y,z)g(\varphi x,u)-
g(\varphi x,z)g(\varphi y,u)-2g(\varphi x,y)g(\varphi z,u)\\
[2mm]
&-g(y,z)\tilde\eta(x)\tilde\eta(u)-g(x,u)\tilde\eta(y)\tilde\eta(z)\\
[2mm]
&+g(x,z)\tilde\eta(y)\tilde\eta(u)+g(y,u)\tilde\eta(x)\tilde\eta(z)],
\quad x,y,z,u \in {\X}\Delta.
\end{array} \leqno (5.9)$$

Comparing (5.9) with the equality from Proposition \ref {P:2.1} we
obtain
$$c+3{\alpha}^2=a+k^2,\quad c-{\alpha}^2=a. \leqno{(5.10)}$$
Hence any integral submanifold $Q^{2n-1}$ of $\Delta$ is an
$\alpha$-Sasakian space form with
$$\alpha = \frac{k}{2}, \quad c=a+\frac{k^2}{4}.$$
Now the relation (5.10) gives the assertion. \hfill{\bf QED}
\vskip 1mm The above statement allows us to obtain examples of
$\alpha$-Sasakian (Sasakian) manifolds of constant
$\varphi$-holomorphic sectional curvatures $c$ satisfying the
condition $c+3\alpha^2<0$ ($c+3<0$) as hypersurfaces of the
K\"ahler space form $(\mathbb{T}^{n-1}_1, g, J, D)$, $g$ given by
$(4.24)$.

Let $({\T}^{n-1}_1,h',J,D)$ be the canonical example of a flat
K\"ahler-Lorentz manifold with $B_0$-distribution $D$ and $g$ be
the K\"ahler metric of constant holomorphic sectional curvatures
$-1$, given by (4.24). We denote by $H^{2(n-1)}_1(O,r)$ any
hypersphere in ${\T}^{n-1}_1$, centered at the origin $O$ and with
radius $r>1$, given by
$$H_1^{2(n-1)}(O,r)=\{\textbf{Z}\in{\T}^{n-1}_1 \,
\vert \; h'(\textbf{Z},\textbf{Z})=-r^2\}.$$ Then an easy
verification shows that $H^{2(n-1)}_1(O,r)$ with the induced from
$({\T}^{n-1}_1,g,J,D)$ structure
$(g,\varphi,\tilde\xi,\tilde\eta)$ is an $\alpha$-Sasakian
manifold with constant $\varphi$-holomorphic sectional curvatures
$c$ such that
$$\alpha=\frac{1}{2r}, \quad c+3\alpha^2=-\frac{r^2-1}{r^2}.$$

Further we give a direct construction of examples of Sasakian
structures with prescribed $\varphi$-holomorphic sectional
curvatures $c$ of type $c+3<0$ using as a base the hypersphere
$H^{2(n-1)}_1(O,r=1)=H^{2(n-1)}_1(1)$.

Let $(h',\varphi,\tilde\xi,\tilde\eta)$ be the induced from
$({\T}^{n-1}_1,h',J,D)$ onto $H^{2(n-1)}_1(1)$ $(-1)$-Sasakian
structure with $h'(\tilde\xi,\tilde\xi)=-1$. We introduce the
following family of Riemannian metrics
$$g=q^2(h'+(1+q^2)\tilde\eta\otimes\tilde\eta), \quad q=const >0
\leqno{(5.11)}$$ on $H^{2(n-1)}_1(1)$. Any of these metrics
generates the corresponding unit vector field $\bar\xi$ and 1-form
$\bar\eta$ determined by
$$\bar\xi=\frac{1}{q^2}\tilde\xi, \quad \bar\eta=-q^2\tilde\eta.$$
In a straightforward way we obtain that
$(H^{2(n-1)}_1(1),g,\varphi,\bar\xi,\bar\eta)$ is a Sasakian
manifold. Further, by direct computations we find that the
Sasakian structure $(g,\varphi,\bar\xi,\bar\eta)$ is of constant
$\varphi$-holomorphic sectional curvatures $c$ satisfying the
relation
$$c+3=-\frac{4}{q^2}.$$

Thus we obtained

{\it Examples of Sasakian space forms with prescribed
$\varphi$-holomorphic sectional curvatures $c$ satisfying the
condition $c+3<0$}:
$$(H^{2(n-1)}_1(1),g,\varphi,\bar\xi,\bar\eta): \quad \bar\xi
=-\frac{c+3}{4}\,\tilde\xi,\quad
\bar\eta=\frac{4}{c+3}\,\tilde\eta,$$
$$g=-\frac{4}{c+3}\left(h'+
\frac{c-1}{c+3}\,\tilde\eta\otimes\tilde\eta\right).$$
\section{K\"ahler structures on rotational hypersurfaces}
In this section we consider three types of rotational
hypersurfaces in spaces with definite or indefinite flat metrics,
which will be endowed with K\"ahler structures of quasi-constant
holomorphic sectional curvatures.

In Subsection $\left \{\begin{array}{l} \vspace{1mm}
6.1,\\
\vspace{1mm}
6.2,\\
\vspace{1mm} 6.3
\end{array} \right .$
we show that any rotational hypersurface of type $\left \{
\begin{array}{l}
\vspace{1mm}
I,\\
\vspace{1mm}
II,\\
\vspace{1mm} III
\end{array} \right .$
carries a K\"ahler structure of quasi-constant holomorphic
sectional curvatures with functions \vskip 2mm $\left \{
\begin{array}{l}
\vspace{1mm}
a+k^2>0, \quad a>0,\\
\vspace{1mm}
a+k^2>0, \quad a<0,\\
\vspace{1mm} a+k^2<0, \quad a<0,
\end{array} \right .$
respectively.

We describe the meridians of those rotational hypersurfaces,
whose K\"ahler metrics are Bochner-K\"ahler
(especially of constant holomorphic sectional curvatures).
\vskip 2mm
\subsection{K\"ahler structures on rotational hypersurfaces of type $I$}
In \cite{GM2} we studied the standard $2n$-dimensional
rotational hypersurfaces $M$ in ${\R}^{2n+1}={\C}^n\times {\R}$
having no common points with the axis of revolution $l={\R}$.
Any such hypersurface $M$ is a one-parameter family of spheres
$S^{2n-1}(s),\,s\in I \subset {\R}$, considered as hyperspheres in
${\C}^n$ with corresponding centers on $l$ and radii $t(s)>0$, $s$
being the natural parameter for the meridian. A rotational
hypersurface $M$ satisfying the conditions
$$t(s)>0, \quad t'(s)>0; \quad s\in I$$
is said to be a {\it rotational hypersurface of type I}.

In \cite{GM2} we have shown that any rotational hypersurface $M$
of type I carries a natural K\"ahler structure $(g,J,\xi)$, which has
the following remarkable property.
\begin{thm}\label{T:6.1} \cite{GM2}
Let $M\, (\dim M= 2n\geq 4)$ be a rotational hypersurface of
type I. Then the K\"ahler structure $(g,J,\xi)$ on $M$ is of
quasi-constant holomorphic sectional curvatures with functions
$$a\geq 0, \quad (a+k^2>0).$$
\end{thm}

The curvature tensor $R$ of the metric $g$ has the form
$$R=a\pi+b\Phi+c\Psi,$$
where
$$a=\frac{4(1-t')}{t^2}, \quad b=8\left(\frac{t'-1}{t^2}-
\frac{t''}{2tt'}\right),\quad
c=\frac{4(1-t')}{t^2}+\frac{5t''}{2tt'}+
\frac{t''^2-t't'''}{2t'^3}.\leqno{(6.1)}$$

In this subsection we describe the rotational hypersurfaces of type I,
whose K\"ahler structure is Bochner flat.

We recall that the Bochner curvature tensor $B(R)$ of a K\"ahler
manifold $(M,g,J)\;(\dim M=2n\geq 4)$ with curvature tensor $R$,
Ricci tensor $\rho$ and scalar curvature $\tau$ is given by
$$\begin{array}{ll}
(B(R))_{\alpha\bar\beta\gamma\bar\delta}=
&R_{\alpha\bar\beta\gamma\bar\delta}- \displaystyle{
\frac{1}{n+2}(g_{\alpha\bar\beta}\rho_{\gamma\bar\delta}
+g_{\gamma\bar\beta}\rho_{\alpha\bar\delta} +
g_{\gamma\bar\delta}\rho_{\alpha\bar\beta}
+ g_{\alpha\bar\delta}\rho_{\gamma\bar\beta})}\\
[3mm]& \displaystyle{+
\frac{\tau}{2(n+1)(n+2)}(g_{\alpha\bar\beta}g_{\gamma\bar\delta}+
g_{\gamma\bar\beta}g_{\alpha\bar\delta})}
\end{array}\leqno{(6.2)}$$
in local holomorphic coordinates.

The manifold $(M,g,J)$ is said to be Bochner flat (or the metric
$g$ is Bochner-K\"ahler) if $B(R)=0$.

\begin{lem}\label{L:6.1}
A K\"ahler manifold whose curvature tensor is of the form
$$R=a\pi+b\Phi+c\Psi,\leqno{(6.3)}$$
is Bochner flat if and only if $c=0$.
\end{lem}
{\it Proof.} Applying the Bochner operator (6.2) to the tensor
(6.3) we find
$$B(R)=c \, \left(\frac{2}{(n+1)(n+2)}\pi - \frac{4}{n+2}\Phi + \Psi\right),$$
which gives the assertion. \hfill{\bf QED} \vskip 1mm

Any rotational hypersurface $M$ of type I is geometrically
determined by the equation $t=t(s)$ \, (or equivalently $s=s(t)$).

\begin{prop}\label{P:6.3}
Let $M$ be a rotational hypersurface of type I. Then the K\"ahler
structure $(g,J)$ is Bochner-K\"ahler if and only if
$$s(t)=\int \frac{dt}{c_1t^4+c_2t^2+1},$$
where $c_1=const, c_2=const$.
\end{prop}

{\it Proof.} Taking into account (6.1) we have
$$c=-\frac{t}{2t'}\left(\frac{t''}{tt'}+4\frac{1-t'}{t^2}\right)'.$$
According to Lemma \ref{L:6.1} we have to solve the equation
$c=0$, i.e.
$$\frac{t''}{tt'}+4\frac{1-t'}{t^2}=const=-2c_2.\leqno{(6.4)}$$
The general solution of (6.4) is
$$s(t)=\int \frac{dt}{c_1t^4+c_2t^2+1}$$
for some constant $c_1$. \hfill{\bf QED}
\vspace{2mm}

We note that the case $c_1=0$ gives the K\"ahler metrics of
constant holomorphic sectional curvatures $a=const>0$ described in
\cite{GM2}.

\subsection{K\"ahler structures on rotational hypersurfaces of type $II$}

Let $({\C}^n,g',J_0)=$ \; \; $({\R}^{2n},g',J_0)$ be the complex
space with the standard complex structure $J_0$ and flat definite
metric $g'$. Further, let $O\textbf{e}$ be a coordinate system on
${\R}$ with the inner product determined by $\textbf{e}^2=-1$ and
$l={\R}$ be the axis of revolution in the space ${\C}^n\times
{\R}$. We denote the product metric in ${\R}^{2n}_1={\C}^n\times
{\R}$ by the same letter $g'$. Then $g'(\textbf{e},\textbf{e})=-1$
and $g'$ is of signature $(2n,1)$.

We consider the class of rotational hypersurfaces having no common
points with the axis of revolution $l$. Then any such hypersurface
$M$ is a one-parameter family of spheres $S^{2n-1}(s),\,s\in I $
considered as hyperspheres in ${\C}^n$ with corresponding centers
$q(s)\textbf{e}$ on $l$ and radii $t(s)>0$. If $\textbf{Z}$ is the
radius vector of any point $p \in M$ with respect to the origin
$O$, then the unit normal $\textbf{n}$ of the parallel
$S^{2n-1}(s)$ at the point $p$ is
$$\textbf{n}=\frac{\textbf{Z}-q(s)\textbf{e}}{t(s)}.$$
Hence
$$\textbf{Z}=t(s)\textbf{n}+q(s)\textbf{e} \leqno{(6.5)}$$
and the meridian $\gamma$ of $M$ is
$$\gamma: \textbf{z}(s)=t(s)\textbf{n} + q(s)\textbf{e} \leqno{(6.6)}$$
in the plane $O\textbf{ne}$ ($\textbf{n}$ - fixed).

Because of (6.6) and (6.5) the tangent vector field $\bar\xi$ to
$\gamma$ is
$$\bar \xi=\frac{d\textbf{z}}{ds}=t'\textbf{n}+q'\textbf{e}=
\frac{\partial \textbf{Z}}{\partial s}. \leqno{(6.7)}$$

We consider rotational hypersurfaces whose meridian $\gamma$ has a
space-like tangent at any point and assume that $s$ is a natural
parameter for $\gamma$, i.e.
$$g'(\frac{d\textbf{z}}{ds},\frac{d\textbf{z}}{ds})=t'^2-q'^2=1.$$

Since the normal to $M$ lies in the plane $O\textbf{ne}$, we
choose the time-like unit vector field $N$ normal to $M$ by the
condition that the couples $(\textbf{n},\textbf{e})$ and
$(\bar\xi,N)$ have the same orientation. Then taking into account
(6.7), we have
$$N=q'\textbf{n}+t'\textbf{e}.$$

\begin{defn}\label{D:6.1}
A rotational hypersurface $M$ in ${\R}^{2n}_1={\C}^n\times{\R}$,
which has no common points with the axis of revolution $l$, is
said to be {\it of type II} if its normals are time-like.
\end{defn}

Let $\nabla'$ be the flat Levi-Civita connection of the metric
$g'$ in ${\R}^{2n}_1={\C}^n\times {\R}$. We denote the induced
definite metric on $M$ by $\bar g$. Let $\bar\eta$ be the 1-form
corresponding to the space-like unit vector field $\bar\xi$ with
respect to the metric $\bar g$, i.e. $\bar\eta(X)=\bar
g(\bar\xi,X), \,\, X\in {\X}M$. If $\bar\nabla$ is the Levi-Civita
connection on $(M,\bar g)$ we have:
$$\begin{array}{l}
\vspace{2mm}
\displaystyle{\nabla'_XY=\bar\nabla_XY+
\left(\frac{\sqrt{t'^2-1}}{t}\bar g(X,Y)+
\frac{1-t'^2+tt''}{t\sqrt{t'^2-1}}\bar\eta(X)\bar\eta(Y)\right)N},
\quad X,Y\in{\X}M;\\
\vspace{2mm}
\bar\nabla_{\bar\xi}\bar\xi=0; \quad
\bar\nabla_x\bar\xi=\displaystyle{\frac{t'}{t}\,x},
\quad \bar g(x,\bar\xi)=0,\quad x\in {\X}M.
\end{array}
\leqno{(6.8)}$$

Then the curvature tensor $\bar R$ of the rotational hypersurface
$(M,\bar g)$ of type II has the form:
$$\bar R=-\frac{t'^2-1}{t^2}\,\bar\pi-\frac{1-t'^2+tt''}{t^2}\,\bar\Phi.
\leqno{(6.9)}$$
This equality implies that the rotational
hypersurface $(M,\bar g)$ of type II is conformally flat. More
precisely, $(M,\bar g, \bar \xi)$ is a subprojective Riemannian
manifold with horizontal sectional curvatures
$\displaystyle{-\frac{t'^2-1}{t^2}}\leq 0$ (cf \cite{GM1}).

As in \cite{GM2}, we consider the almost contact Riemannian
structure $(\varphi, \tilde{\bar\xi}, \tilde{\bar\eta},\bar g)$
on the parallels $S^{2n-1}(s),\;s\in I$ of the
rotational hypersurface $M$ and obtain that any parallel is
$\displaystyle{\frac{1}{t}}$-Sasakian.

This allows us to introduce the almost complex structure $J$ on
$(M,\bar g)$ subordinated to the orientation $\bar\xi$ of the
meridians by
$$J_{|D}:=J_0, \quad J\bar\xi:=\tilde{\bar\xi}, \quad J\tilde{\bar\xi}
:=-\bar\xi.\leqno{(6.10)} $$

Similarly to the definite case \cite{GM2} we have

\begin{prop}\label{P:6.4}
Let $(M,\bar g)$ be a rotational hypersurface of type II in
${\R}^{2n}_1={\C}^n\times \mathbb{R}$ whose meridians are oriented
with the space-like unit vector field $\bar\xi$ . If $J$ is the
almost complex structure $(6.10)$ associated with $\bar\xi$, then
the covariant derivative of $J$ satisfies the identity
$$(\bar\nabla_XJ)Y=\frac{t'-1}{t}\left(\bar g(X,Y)\tilde{\bar\xi}-
\tilde{\bar\eta}(Y)X-\bar\eta(Y)JX+\bar g(JX,Y)\bar\xi\right)
\leqno{(6.11)}$$
for all vector fields $X,Y\in{\X}M$.
\end{prop}

The identity (6.11) shows that $(M,\bar g,J)$ is a locally
conformal K\"ahler manifold in all dimensions $2n\geq 4$ with Lee
form $\displaystyle{\frac{1-t'}{t}\,\bar\eta}$.

Our aim in this subsection is to define a nontrivial K\"ahler
metric on $(M,\bar g,J)$, which is naturally determined by its
geometric structures.

If $(M,\bar g,J)$ is a rotational hypersurface of type II, then
$t'^2\geq 1$. Therefore we can always choose the orientation
$\bar\xi$ of the meridians so that $t'\geq 1$.

In what follows we assume that
$$t(s)>0,\quad t'(s)\geq 1; \quad s\in I. \leqno{(6.12)}$$

Under the conditions (6.12) we construct the structure $(g, \xi)$:
$$g=\bar g+(t'-1)(\bar\eta\otimes\bar\eta+\tilde{\bar\eta}
\otimes\tilde{\bar\eta}), \quad \xi=\frac{1}{\sqrt{t'}}\,\bar\xi,
\quad \eta=\sqrt{t'}\,\bar\eta. \leqno{(6.13)}$$

Taking into account (6.11) we obtain that the K\"ahler form of the
metric (6.13) is closed, i.e. $g$ is a K\"ahler metric. More
precisely, we have
\begin{thm}\label{T:6.2}
Let $(M,\bar g,J,\bar\xi)\,(2n\geq 4)$ be a rotational
hypersurface of type II and assume that $(6.12)$ hold good. Then
the  K\"ahler metric $g$, given by $(6.13)$, is of quasi-constant
holomorphic sectional curvatures with functions
$$a\leq 0, \qquad a+k^2>0.$$
\end{thm}

{\it Proof.} Calculating the relation between the connections of
the metrics in (6.13) in view of (6.9) we find the curvature
tensor $R$ of the K\"ahler metric $g$:
$$R=a\pi+b\Phi+c\Psi,$$
where
$$a=\frac{4(1-t')}{t^2}, \quad b=8\left(\frac{t'-1}{t^2}-
\frac{t''}{2tt'}\right), \quad
c=\frac{4(1-t')}{t^2}+\frac{5t''}{2tt'}+\frac{t''^2-t't'''}{2t'^3}.
\leqno{(6.14)}$$
Applying Proposition 2.3 \cite{GM2} we obtain
that $(M,g,J,\xi)$ is of quasi-constant holomorphic sectional
curvatures.

Since $t'\geq 1$, then we have $a \leq 0$.

From (6.8) and the relation between the connections of $g$ and
$\bar g$  it follows that
$$\displaystyle{\nabla_x\xi=
\frac{\sqrt{t'}}{t}\,x-\frac{t''}{2t'\sqrt{t'}}\,\eta(Jx)J\xi}$$
for all $x\in{\X}M,$ $g(\xi,x)=0$. According to (2.1) the function
$k$ of the structure $(g,J,\xi)$ is
$\displaystyle{k=2\frac{\sqrt{t'}}{t}}$. Taking into account
(6.14) we find
$$a+k^2=\frac{4}{t^2}>0.$$
\hfill{\bf QED}
\vskip 1mm
As a consequence of Theorem \ref{T:6.2}
we can find the rotational hypersurfaces $(M,\bar g, J)$ of type
II whose K\"ahler metric (6.13) is of constant holomorphic
sectional curvatures.

Let $b=0$ in (6.14). Then Corollary 3.6 \cite{GM2} implies that
$c=0$ and the metric $g$ is of constant holomorphic sectional
curvatures $a=const \leq 0$.

Solving the equation
$$b=8\left(\frac{t'-1}{t^2}-\frac{t''}{2tt'}\right)=0$$
we obtain the meridian in the form $q=q(t)$.

Namely, we have
\begin{prop}\label{P:6.5}
Any rotational hypersurface $(M,\bar g,J)$ of type II, whose
K\"ahler metric $(6.13)$ is of constant holomorphic sectional
curvatures $a=const <0$ is generated by a meridian of the type
$$\gamma : q=\pm\frac{1}{\sqrt{-a}}\left(\sqrt{8-at^2}+
\ln{\frac{\sqrt{8-at^2}- 2}{\sqrt{8-at^2}+2}}\right)+q_0, \quad t
> 0$$ in the hyperbolic plane $O\textbf{ne}$.
\end{prop}

Similarly to Proposition \ref{P:6.3} we obtain the following
statement.
\begin{prop}\label{P:6.6}
Let $(M,\bar g,J)$ be a rotational hypersurface of type II
generated by the meridian
$$\gamma: \textbf{z}(s)=t(s)\textbf{n}+q(s)\textbf{e}, \quad s\in I$$
in the hyperbolic plane $O\textbf{ne}$. Then the metric $g$ given
by $(6.13)$ is Bochner-K\"ahler if and only if
$$s(t)=\int \frac{dt}{c_1t^4+c_2t^2+1},$$
where $c_1=const, c_2=const$.
\end{prop}

We note that the case $c_1=0$ is described in Proposition
\ref{P:6.5}.

\subsection{K\"ahler structures on rotational hypersurfaces of type $III$}

Let $({\C}^n,h',J_0)=$ \; \; $({\R}^{2(n-1)}_2,h',J_0)$ be the
K\"ahler-Lorentz space with the standard complex structure $J_0$
and flat indefinite metric $h'$ of signature $(2(n-1),2)$.
Further, let $O\textbf{e}$ be a coordinate system on ${\R}$ with
the inner product determined by $\textbf{e}^2=+1$ and $l={\R}$ be
the axis of revolution in the space ${\R}^{2(n-1)}_2\times {\R}
={\C}^n\times {\R}$. We denote the product metric in
${\R}^{2n-1}_2={\C}^n\times {\R}$ by the same letter $h'$. Then
$h'(\textbf{e},\textbf{e})=+1$ and $h'$ is of signature
$(2(n-1),2)$.

We consider rotational hypersurfaces $M$ with parallels
$H^{2(n-1)}_1$, which are hyperspheres with respect to the metric
$h'$ in the time-like domain ${\T}^{n-1}_1\subset {\C}^n$. Then
$M$ is a one-parameter family of spheres $H^{2(n-1)}_1(s),\,s\in I
$ with corresponding centers $q(s)\textbf{e}$ on $l$ and radii
$t(s)>0$. If $\textbf{Z}$ is the radius vector of any point $p \in
M$ with respect to the origin $O$, then the unit normal
$\textbf{n}$ of the parallel $H^{2(n-1)}_1(s)$ at the point $p$ is
$$\textbf{n}=\frac{\textbf{Z}-q(s)\textbf{e}}{t(s)},
\quad h'(\textbf{n},\textbf{n})=-1.$$ Hence
$$\textbf{Z}=t(s)\textbf{n}+q(s)\textbf{e} \leqno{(6.15)}$$
and the meridian $\gamma$ of $M$ is
$$\gamma: \textbf{z}(s)=t(s)\textbf{n} + q(s)\textbf{e} \leqno{(6.16)}$$
in the plane $O\textbf{ne}$ ($\textbf{n}$ - fixed).

Because of (6.16) and (6.15) the tangent vector field $\bar\xi$ to
$\gamma$ is
$$\bar \xi=\frac{d\textbf{z}}{ds}=t'\textbf{n}+q'\textbf{e}
=\frac{\partial \textbf{Z}}{\partial s}. \leqno{(6.17)}$$

We consider rotational hypersurfaces whose meridian $\gamma$ has a
time-like tangent at any point and assume that $s$ is a natural
parameter for $\gamma$, i.e.
$$h'(\frac{d\textbf{z}}{ds},\frac{d\textbf{z}}{ds})=-t'^2+q'^2=-1.$$

Since the normal to $M$ lies in the plane $O\textbf{ne}$, we
choose the space-like unit vector field $N$ normal to $M$ by the
condition that the couples $(\textbf{n},\textbf{e})$ and
$(\bar\xi,N)$ have the same orientation. Then taking into account
(6.17), we have
$$N=q'\textbf{n}+t'\textbf{e}.$$

\begin{defn}\label{D:6.2}
A rotational hypersurface $M$ in ${\R}^{2n-1}_2={\C}^n\times{\R}$,
which has no common points with the axis of revolution $l={\R}$,
is said to be {\it of type III} if its normals are space-like.
\end{defn}

Let $\nabla'$ be the flat Levi-Civita connection of the metric
$h'$ in ${\R}^{2n-1}_2={\C}^n\times {\R}$. We denote by $\bar h$
the induced indefinite metric on $M$ of signature $(2(n-1),2)$.
Let $\bar\eta$ be the 1-form corresponding to the unit time-like
vector field $\bar\xi$ with respect to the metric $\bar h$, i.e.
$\bar\eta(X)=\bar h(\bar\xi,X), \,\, X\in {\X}M$. If $\bar\nabla$
is the Levi-Civita connection on $(M,\bar h)$ we have:
$$\begin{array}{l}
\vspace{2mm}
\displaystyle{\nabla'_XY=\bar\nabla_XY-
\left(\frac{\sqrt{t'^2-1}}{t}\,\bar h(X,Y)+
\frac{-1+t'^2+tt''}{t\sqrt{t'^2-1}}\bar\eta(X)\bar\eta(Y)\right)N},
\; X,Y\in{\X}M;\\
\vspace{2mm}
\bar\nabla_{\bar\xi}\bar\xi=0; \quad
\bar\nabla_x\bar\xi=\displaystyle{\frac{t'}{t}\,x},
\quad \bar h(x,\bar\xi)=0, \quad x\in {\X}M.
\end{array}
\leqno{(6.18)}$$

Then the curvature tensor $\bar R$ of the rotational hypersurface
$(M,\bar h)$ of type III has the form:
$$\bar R=\frac{t'^2-1}{t^2}\,\bar\pi+\frac{-1+t'^2+tt''}{t^2}\,\bar\Phi,
\leqno{(6.19)}$$
where $\bar\pi$ and $\bar\Phi$ are the tensors
$$\begin{array}{ll}
\vspace{2mm}
\bar\pi (X,Y)Z= & \bar h(Y,Z)X-\bar h(X,Z)Y,\\
\vspace{2mm}
\bar \Phi (X,Y)Z= & \bar h(Y,Z)\bar\eta(X)\bar\xi-
\bar h(X,Z)\bar\eta(Y)\bar\xi\\
\vspace{1mm}
& +\bar\eta(Y)\bar\eta(Z)X-\bar\eta(X)\bar\eta(Z)Y,
\quad X,Y,Z \in {\X}M.
\end{array}$$
The equality (6.19) implies that the rotational hypersurface
$(M,\bar h)$ of type III is conformally flat.

Now we consider the almost contact Riemannian structure $(\varphi,
\tilde{\bar\xi}, \tilde{\bar\eta},\bar h)$ on the parallel
$H^{2(n-1)}_1(s),\;s\in I$ of the rotational hypersurface $(M,\bar
h)$ which arises in a similar way as in the definite case (cf
\cite{T1,T2}):

$$\begin{array}{l}
\vspace{2mm}
\tilde{\bar\xi}:=J_0n, \qquad \tilde{\bar\eta}(x)
:=\bar h(\tilde{\bar\xi},x);\\
\vspace{2mm}
\varphi x:= J_0x-\tilde{\bar\eta}(x)n, \quad x\in
{\X}H^{2(n-1)}_1(s).
\end{array}\leqno{(6.20)}$$

It is clear that $\bar h(\tilde{\bar\xi},\tilde{\bar\xi})=-1$. The
relations (6.20) imply that
$$\begin{array}{l}
\vspace{2mm}
\varphi \, \xi =0;\quad
\varphi^2 x=-x-\tilde{\bar\eta}(x)\tilde{\bar\xi};\\
\vspace{2mm}
\bar h(\varphi x,\varphi y)=\bar h(x,y)
+\tilde{\bar\eta}(x)\tilde{\bar\eta}(y),
\quad x,y \in {\X}H^{2(n-1)}_1(s).
\end{array}$$

Let us denote by $\mathcal{D}$ the induced Levi-Civita connection
of the metric $\bar h$ on $H^{2(n-1)}_1(s)$ as a submanifold of
${\T}^{n-1}_1(s)\subset {\C}^n$. Then the Weingarten and Gauss
formulas of the imbedding $H^{2(n-1)}_1(s)\subset {\C}^n$ are:
$$\begin{array}{l}
\vspace{2mm}
\displaystyle{\nabla'_x\textbf{n}=\frac{1}{t}\,x;}\\
\vspace{2mm}
\displaystyle{\nabla'_xy=\mathcal{D}_xy+\frac{1}{t}\,\bar
h(x,y)\textbf{n}, \quad x,y\in
{\X}H^{2(n-1)}_1(s).}
\end{array}\leqno{(6.21)}$$

From (6.20) and (6.21) we obtain consequently
$$\begin{array}{l}
\vspace{2mm}
\displaystyle{\nabla'_x\tilde{\bar\xi}=
\frac{1}{t}(\varphi x+\tilde{\bar\eta}(x)\textbf{n});}\\
\vspace{2mm}
\displaystyle{\mathcal{D}_x\tilde{\bar\xi}=\frac{1}{t}\varphi x,}
\quad x\in {\X}H^{2(n-1)}_1(s).\\
\end{array}\leqno{(6.22)}$$

Let $T_pM$ be the tangent space to $M$ at any point $p\in M$. Then
the vector fields $\bar\xi$ and $\tilde{\bar\xi}$ defined by
(6.20) determine a distribution $D$ such that $D^{\perp}=
span\{\bar\xi,\tilde{\bar\xi}\}$. The distribution $D$ is
space-like, while the distribution $D^{\perp}$ is time-like.

We define an almost complex structure $J$ on $(M,\bar h)$
subordinated to the orientation $\bar\xi$ of the meridians
$\gamma$ as follows:
$$J_{|D}:=J_0, \quad J\bar\xi:=\tilde{\bar\xi}, \quad J\tilde{\bar\xi}
:=-\bar\xi. \leqno{(6.23)}$$

Similarly to Proposition \ref{P:6.4} we have

\begin{prop}\label{P:6.7}
Let $(M,\bar h)$ be a rotational hypersurface of type III in
${\R}^{2n-1}_2={\C}^n\times \mathbb{R}$ whose meridians $\gamma$
are oriented with the time-like unit vector field $\bar\xi$ . If
$J$ is the almost complex structure $(6.23)$ associated with
$\bar\xi$, then the covariant derivative of $J$ satisfies the
identity
$$(\bar\nabla_XJ)Y=\frac{1-t'}{t}\left(\bar h(X,Y)\tilde{\bar\xi}-
\tilde{\bar\eta}(Y)X-\bar\eta(Y)JX+\bar h(JX,Y)\bar\xi\right)
\leqno{(6.24)}$$
for all vector fields $X,Y\in{\X}M$.
\end{prop}

{\it Proof.} We calculate the components of $(\bar\nabla_XJ)Y, \;
X,Y \in {\X}M$:
$$\begin{array}{l}
\vspace{2mm}
\displaystyle{(\bar\nabla_xJ)y=\frac{1-t'}{t}(\bar
h(x,y)\tilde{\bar\xi}
+\bar h(\varphi x,y)\bar\xi-\tilde{\bar\eta}(y)n)},\\
\vspace{2mm}
\displaystyle{(\bar\nabla_xJ)\bar\xi=\frac{1-t'}{t}\varphi x,}
\quad
x,y \in {\X}M, \quad \bar h(\bar\xi,x)=\bar h(\bar\xi,y)=0;\\
\vspace{2mm}
(\bar\nabla_{\bar\xi}J)x_0=0, \quad x_0 \in {\X}M,
\quad \bar h(x_0,\bar\xi)=\bar h(x_0,\tilde{\bar\xi})=0;\\
\vspace{2mm}
(\bar\nabla_{\bar\xi}J)\tilde{\bar\xi}=0, \quad
(\bar\nabla_{\bar\xi}J)\bar\xi=0.\\
\end{array}$$

These equalities imply the assertion. \hfill{\bf QED} \vspace{2mm}

The identity (6.24) shows that $(M,\bar h,J)$ is a locally
conformal K\"ahler manifold in all dimensions $2n\geq 4$ with Lee
form $\displaystyle{\frac{1-t'}{t}\,\bar\eta}$. This implies that
$(M,\bar h,J)$ carries a conformal K\"ahler metric of signature
$(2(n-1),2)$ which is flat.

Our aim in this subsection is to define a nontrivial definite
K\"ahler metric $g$ on $(M,\bar h,J)$, which is naturally
determined by its geometric structures.

If $(M,\bar h,J)$ is a rotational hypersurface of type III, then
$t'^2\geq 1$. Therefore we can always choose the orientation
$\bar\xi$ of the meridians so that $t'\leq -1$.

In what follows we assume that
$$t(s)>0,\quad t'(s)\leq -1; \quad s\in I. \leqno{(6.25)}$$

Under the conditions (6.25) we construct the structure $(g,\xi)$:
$$g=\bar h+(1-t')(\bar\eta\otimes\bar\eta+\tilde{\bar\eta}
\otimes\tilde{\bar\eta}), \quad \xi=\frac{1}{\sqrt{-t'}}\,\bar\xi,
\quad \eta=-\sqrt{-t'}\,\bar\eta. \leqno{(6.26)}$$

Taking into account the defining condition (6.26) and (6.25), we
obtain that $g$ is a definite metric and $\xi$ is a unit vector
field. Because of (6.26) and (6.24) it follows that $g$ is a
K\"ahler metric on $M$.

More precisely, we have
\begin{thm}\label{T:6.3}
Let $(M,\bar h,J,\bar\xi)\,(2n\geq 4)$ be a rotational
hypersurface of type III and assume that $(6.25)$ hold good. Then
the  K\"ahler metric $g$, given by $(6.26)$, is of quasi-constant
holomorphic sectional curvatures with functions
$$a<0, \qquad a+k^2<0.$$
\end{thm}

{\it Proof.} Let $\nabla$ be the Levi-Civita connection of the
metric (6.26). We calculate the relation between $\bar\nabla$ and
$\nabla$:
$$\begin{array}{ll}
\vspace{2mm}
\bar\nabla_XY = & \nabla_XY -
\displaystyle{\frac{t''}{2t'\sqrt{-t'}}
\{[\eta(X)\eta(Y)-\eta(JX)\eta(JY)]\xi
-[\eta(X)\eta(JY)+\eta(JX)\eta(Y)]J\xi\}}\\
\vspace{2mm}
& \displaystyle{-\frac{1-t'}{t\sqrt{-t'}}
\{\eta(JX)JY+\eta(JY)JX-[\eta(X)\eta(JY)+
\eta(JX)\eta(Y)]J\xi}\\
\vspace{2mm}
&-t'[g(X,Y)-\eta(JX)\eta(JY)-\eta(X)\eta(Y)]\xi-2\eta(JX)\eta(JY)\xi\}
\end{array}$$
for all $X,Y \in {\X}M$.

Taking into account (6.18) we find
$$\nabla_X\xi=\displaystyle{\frac{t'}{t\sqrt{-t'}}(X-\eta(X)\xi)-
\frac{t''}{2t'\sqrt{-t'}}\eta(JX)J\xi}, \quad X \in {\X}M.
\leqno{(6.27)}$$

Then we find the curvature tensor $R$ of the K\"ahler metric $g$:
$$R=a\pi+b\Phi+c\Psi,$$
where
$$a=\frac{4(t'-1)}{t^2}, \quad b=-8\left(\frac{t'-1}{t^2}-
\frac{t''}{2tt'}\right), \quad
c=\frac{4(t'-1)}{t^2}-2\frac{t''}{tt'}-\frac{t''^2}{2tt'^3}
\left(\frac{tt'}{t''}\right)'. \leqno{(6.28)}$$

Applying Proposition 2.3 \cite{GM2} we obtain that $(M,g,J,\xi)$
is of quasi-constant holomorphic sectional curvatures.

Since $t'\leq -1$, then we have $a < 0$.

From (6.27) it follows that the function $k$ of the structure
$(g,J,\xi)$ is $\displaystyle{k=\frac{2t'}{t\sqrt{-t'}}}$. Taking
into account (6.28), we find
$$a+k^2=-\frac{4}{t^2}<0.$$
\hfill{\bf QED} \vskip 2mm

As a consequence of Theorem \ref{T:6.3} we can find the rotational
hypersurfaces $M$ of type III whose K\"ahler metric (6.26) is of
constant holomorphic sectional curvatures.

Let $b=0$ in (6.28). Then Corollary 3.6 \cite{GM2} implies that
$c=0$ and the metric $g$ is of constant holomorphic sectional
curvatures $a=const < 0$.

Solving the equation
$$b=-8\left(\frac{t'-1}{t^2}-\frac{t''}{2tt'}\right)=0,$$
we obtain the meridian in the form $q=q(t)$.

Namely, we have
\begin{prop}\label{P:6.7}
Any rotational hypersurface $(M,\bar h,J)$ of type III, whose
K\"ahler metric $(6.28)$ is of constant holomorphic sectional
curvatures $a=const<0$, is generated by a meridian of the type
$$\gamma : q=\frac{1}{-a}\left(\sqrt{a(8+at^2)}-
2\sqrt{-a}\, \arctan{\frac{1}{2}\sqrt{-(8+at^2)}}\right), \quad t
> \frac{2\sqrt{2}}{\sqrt{-a}}$$ in the hyperbolic plane
$O\textbf{ne}$.
\end{prop}

Similarly to Proposition \ref{P:6.3} and Proposition \ref {P:6.6}
we obtain the following statement.
\begin{prop}\label{P:6.8}
Let $(M,\bar h,J)$ be a rotational hypersurface of type III
generated by the meridian
$$\gamma: \textbf{z}(s)=t(s)\textbf{n}+q(s)\textbf{e}, \quad s\in I$$
in the hyperbolic plane $O\textbf{ne}$. Then the metric $g$ given
by $(6.28)$ is Bochner-K\"ahler if and only if
$$s(t)=\int \frac{dt}{c_1t^4+c_2t^2+1},$$
where $c_1=const, c_2=const$.
\end{prop}

We note that the case $c_1=0$ is described in Proposition
\ref{P:6.7}.


\begin{thebibliography}{99}
\bibitem{B}
Bryant, R. {\it Bochner-K\"ahler metrics,} J. Amer. Math. Soc.,
{\bf 14} (2001), 623-715.
\bibitem{GM1}
Ganchev, G.; Mihova, V. {\it Riemannian manifolds of
quasi-constant sectional curvature}, J. reine und angew. Math.,
{\bf 522} (2000), 119-141.
\bibitem{GM2}
Ganchev, G.; Mihova, V. {\it K\"ahler manifolds of quasi-constant
holomorphic sectional curvatures}, ArXiv: math.DG/0505671, to
appear.
\bibitem{JV}
Janssens, D.; Vanhecke, L. {\it Almost contact structures and
curvature tensors}, Kodai Math. J., {\bf 4} (1981), 1-27.
\bibitem{O}
Ogiue, K. {\it On almost contact manifolds admitting axiom of
planes or axiom of free mobility,} Kodai Math. Sem. Rep., {\bf 16}
(1964), 223-232.
\bibitem{TL}
Tachibana, S.; Liu, R.C. {\it Notes on K\"ahlerian metrics with
vanishing Bochner curvature tensor}, Kodai Math. Sem. Rep., {\bf
22} (1970), 313-321.
\bibitem{T}
Tanno, S. {\it Sasakian manifolds with constant
$\varphi$-holomorphic sectional curvature,} T\^ohoku Math. J.,
{\bf 21} (1969), 501-507.
\bibitem{T1}
Tashiro, Y. {\it On contact structures on hypersurfaces in almost
complex manifolds I,} T\^ohoku Math. J., {\bf 15} (1963), 62-79.
\bibitem{T2}
Tashiro, Y. {\it On contact structures on hypersurfaces in almost
complex manifolds II,} T\^ohoku Math. J., {\bf 15} (1963),
167-175.
\bibitem{TV}
Tricerri, F.; Vanhecke, L. {\it Curvature tensors on almost
Hermitian manifolds}, Trans. Amer. Math. Soc., {\bf 267} (1981),
365-398.


\end{thebibliography}
\end{document}